\documentclass[11pt]{article}

\usepackage{graphicx}
\usepackage{latexsym,amsmath,amsfonts,amscd, amsthm, dsfont}
\usepackage{bm,color}
\usepackage{epsfig,verbatim,epstopdf,graphics}
\usepackage{subfigure}
\usepackage{changebar}
\usepackage{multirow}
\usepackage{hyperref}

\usepackage{algorithmic}% http://ctan.org/pkg/algorithms
\usepackage{yhmath}%�ӻ���ͷ��������  \wideparen
\usepackage{booktabs} %\cmidrule(lr){2-5} \cmidrule(lr){6-9}
\usepackage{tikz}
\usepackage{verbatim}
\usetikzlibrary{arrows,backgrounds,snakes,shapes}
\numberwithin{equation}{section}

\graphicspath{{./}{./figures/}}
\allowdisplaybreaks

\newtheoremstyle{plainNoItalics}{}{}{\normalfont}{}{\bfseries}{.}{ }{}

\theoremstyle{plain}
\newtheorem{thm}{Theorem}[section]

\theoremstyle{plainNoItalics}

\newtheorem{prop}[thm]{Proposition}
\newtheorem{exa}[thm]{Example}

\newcommand{\bx}{{\bf x}}

\newcommand{\bv}{{\bf v}}
\newcommand{\bb}{{\bf b}}

\newcommand{\bl}{{\bf l}}

\newcommand{\ba}{{\bf a}}

\newcommand{\bV}{{\bf V}}

\newcommand{\bj}{{\bf j}}
\newcommand{\bi}{{\bf i}}

\newcommand{\bW}{{\bf W}}
\newcommand{\bk}{{\bf k}}

\newcommand{\bzero}{\mathbf{0}}
\newcommand{\bone}{\mathbf{1}}
\newcommand{\be}{\begin{eqnarray}}
\newcommand{\ee}{\end{eqnarray}}
\newcommand{\beno}{\begin{eqnarray*}}
	\newcommand{\eeno}{\end{eqnarray*}}

\makeatletter

\newcommand{\Rmnum}[1]{\expandafter\@slowromancap\romannumeral #1@}

\newcommand{\abs}[1]{\left\vert#1\right\vert}
\newcommand{\norm}[1]{\left\Vert#1\right\Vert}
\makeatother

\begin{document}
\baselineskip=1.5pc

\vspace{.5in}

\begin{center}

{\large\bf An adaptive sparse grid local discontinuous Galerkin method for Hamilton-Jacobi equations in high dimensions}

\end{center}

\vspace{.1in}

\centerline{
Wei Guo  \footnote{Department of Mathematics and Statistics, Texas Tech University, Lubbock, TX, 70409, USA. E-mail: weimath.guo@ttu.edu. Research is supported by NSF grant DMS-1830838} \qquad
Juntao Huang \footnote{Department of Mathematics, Michigan State University, East Lansing, MI 48824, USA. E-mail: huangj75@msu.edu. Corresponding author} \qquad
Zhanjing Tao  \footnote{School of Mathematics, Jilin University, Changchun, Jilin 130012, China. E-mail: zjtao@jlu.edu.cn}  \qquad
Yingda Cheng  \footnote{Department of Mathematics, Department of Computational Mathematics, Science and Engineering, Michigan State University, East Lansing, MI 48824, USA. E-mail: ycheng@msu.edu. Research is supported by NSF grants DMS-1453661 and DMS-1720023}
}

\vspace{.4in}

\centerline{\bf Abstract}

We are interested in numerically solving the Hamilton-Jacobi (HJ) equations, which arise in optimal control and many other applications. Oftentimes, such equations are posed in high dimensions, and this poses great numerical challenges. This work proposes a class of  adaptive sparse grid (also called adaptive multiresolution) local discontinuous Galerkin (DG) methods for solving Hamilton-Jacobi equations in high dimensions. By using the sparse grid techniques, we can treat moderately high dimensional cases. Adaptivity is   incorporated  to capture kinks and other local structures of the solutions. Two classes of multiwavelets are used to achieve multiresolution,  which are the orthonormal Alpert's multiwavelets and the interpolatory multiwavelets. Numerical tests in up to four dimensions are provided to validate the performance of the method.

\bigskip

\bigskip
%\vfill


%{\bf AMS(MOS) subject classification:} 65M99

\pagenumbering{arabic}

%section 1
\section{Introduction}
\label{sec1}
\setcounter{equation}{0}
\setcounter{figure}{0}
\setcounter{table}{0}

In this paper, we consider the Hamilton-Jacobi (HJ) equation 
\begin{equation}\label{eq:hj}
	\phi_t + H(\nabla \phi) = 0,
\end{equation}
on the bounded domain $[0,1]^d$ in arbitrary $d$ dimension, subject to initial condition $\phi(\bx,0)=\phi_0(\bx)$ and appropriate boundary conditions. The HJ equation has diverse applications in science and engineering,  such as optimal control, seismic waves, crystal growth, robotic navigation, image processing, calculus of variations, among others \cite{lions1982generalized}. In this paper, we develop a class of adaptive sparse grid (also called adaptive multiresolution) discontinuous Galerkin (DG) methods for approximating the viscosity solution of \eqref{eq:hj}. 
The concept of viscosity solution  was developed by Crandal and Lions in \cite{crandall1983viscosity, crandall1984some} to single out the physically relevant weak solution. Under certain assumptions, the viscosity solution can be interpreted by the Hopf formula \cite{evans10}, and the numerical approximation to the viscosity solution is of practical interest \cite{shu2007high}.  It is well known that the viscosity solution of the HJ equation is only Lipschitz continuous 
and may develop discontinuous derivatives in finite time regardless of smoothness of the initial condition. 
Various numerical methods for solving \eqref{eq:hj} have been developed in the literature \cite{shu2007high}, such as
the monotone methods \cite{crandall1984two, abgrall1996numerical, lafon1996high}, 
the essentially non-oscillatory (ENO) methods \cite{osher1988fronts, osher1991high}, the weighted ENO (WENO)   %\cite{}%,ha2017sixth}, 
and the Hermite WENO (HWENO) methods \cite{jiang2000weighted, zhang2003high,qiu2005hermite} among many others. %, qiu2007hermite, zhu2013hermite, tao2017dimension, zheng2016directly, zheng2017finite}. See the review article for more detailed discussion.
In this work, we choose to use
the DG discretization  due to their distinguished advantages in handling geometry, boundary conditions and accommodating adaptivity, which are highly desirable for efficiently solving  the HJ equation. Several DG schemes have been proposed in the literature \cite{hu1999discontinuous, lepsky2000analysis, cheng2007discontinuous, li2010central, yan2011local, guo2011local, cheng2014hj,ke2019alternative}. Here, we use the   local DG (LDG) method developed by Yan and Osher \cite{yan2011local}, which has provable property for the piecewise constant case and is easy to implement, although we remark that the extensions to other DG formulations are also possible. 

Beyond the need to capture the viscosity solution, another major numerical challenge for HJ equation is that it is often posed in high dimensions, and any standard numerical discretization becomes inefficient due to the curse of dimensionality. Recent years have seen a surge of interests in designing numerical  solutions of HJ equations in high dimensions. Various approaches have been proposed, including those using sparse grid \cite{bokanowski2013adaptive,garcke2017suboptimal,kang2017mitigating}, model order reduction \cite{kunisch2004hjb}, tensor decomposition \cite{dolgov2019tensor}, Hopf formula \cite{darbon2016algorithms,chow2017algorithm,chow2019algorithm} and   machine learning \cite{han2018solving,nakamura2019adaptive,darbon2019overcoming,darbon2020some}, to name a few. Some of the work above is feasible for HJ equations in hundreds of dimensions  for some special cases, and continued efforts  to develop efficient numerical solvers for high-dimensional HJ PDEs constitutes a vibrant research area due to their wide applications in control and differential games.

%Hence, adaptivity in space is crucial to capture the viscosity solution efficiently.
In this paper, we take the sparse grid approach \cite{bungartz2004sparse}, which has been used in \cite{bokanowski2013adaptive,garcke2017suboptimal,kang2017mitigating} for   computations in moderately high dimension.
%Based on the DG discretization, this paper considers adaptive HJ simulations in the multiresolution sense. 
The scheme we proposed relies on multiresolution analysis (MRA) \cite{mallat1999wavelet} and is designed to be high order accurate. %explores mesh hierarchy and induces nested polynomial approximation spaces to accelerate the computation. MRA also plays a central role in the popular sparse grid methods \cite{bungartz2004sparse}, which are  effective for simulating high-dimensional problems. We mention a relevant work by Bokanowski, \emph{et al.} \cite{bokanowski2013adaptive} for an adaptive sparse grid semi-Lagrangian method for the Hamilton-Jacobi-Bellman (HJB) equation (a special type of the HJ equation). 
In a line of research, we have developed  a family of adaptive sparse grid (or adaptive multiresolution) DG methods for  linear transport equations with application to kinetic equations \cite{guo2017adaptive}, hyperbolic conservation laws \cite{huang2019adaptive}, and wave equations \cite{huang2020adaptive}. By incorporating MRA and the sparse grid ideas, our methods are able to efficiently capture small-scale structures, and more importantly, work very well in high dimensions. %Hence, the adaptive multiresolution DG framework is ideal for simulating the HJ equation. There are two key ingredients in the proposed  methodology. First, the local DG (LDG) method developed by Yan and Osher \cite{yan2011local} is employed as the base weak formulation, which is proven effective in capturing the correct viscosity solution. Second, 
In particular, in  \cite{huang2019adaptive,huang2020adaptive}, we use two classes of multiwavelets to achieve MRA. The tensor-product Alpert's multiwavelets are used as the DG function space, following the approach developed in \cite{wang2016elliptic,guo2016sparse,guo2017adaptive} for linear equations. Besides, the interpolatory multiwavelets  for MRA quadrature \cite{tao2019collocation} are used for computations of nonlinear terms. %Noteworthy, the LDG formulation is the most amenable to the underlying MRA DG framework in implementation, while other DG formulations can also be used.
 Numerical experiments for benchmark tests in up to four dimension verify the efficiency and efficacy of the method in capturing the viscosity solution of the HJ equations.

The rest of the paper is organized as follows. In Section \ref{sec:mra}, we review the fundamentals of Alpert's and interpolatory
multiwavelets. In Section \ref{sec:dg}, we present the LDG method for the HJ equation with MRA. Some theoretical
results and implementation details are discussed. Section \ref{sec:numer} contains numerical examples. In Section \ref{sec:conclusion}, we include the conclusion of this paper.

%More details about the development of numerical schemes, especially the high order ones, for solving the HJ equation can be found in the review \cite{shu2007high} and the references therein.

% section 2

\section{Multiresolution Analysis and Multiwavelets}
\label{sec:mra}

In this section, we first review the fundamentals of MRA of DG approximation spaces and the associated multiwavelets. Two classes of multiwavelets, namely  the $L^2$ orthonormal Alpert's multiwavelets \cite{alpert1993class} and the interpolatory multiwavelets \cite{tao2019collocation}, are considered. We also introduce a set of key notations used throughout the paper.

\subsection{Alpert's multiwavelets}\label{subsec:alpert-basis}

We start with the construction of Alpert's multiwavelets \cite{alpert1993class}, which have been employed to develop a class of sparse grid DG methods for solving PDEs in high dimensions \cite{wang2016elliptic,guo2016sparse}.
For a unit sized interval $\Omega=[0,1]$,  we define a set of nested grids $\Omega_0,\,\Omega_1,\ldots$, for which the $n$-th level grid $\Omega_n$ consists of $2^n$ uniform cells
\begin{equation*}
I_{n}^j=(2^{-n}j, 2^{-n}(j+1)], \quad j=0, \ldots, 2^n-1,
\end{equation*}
 Denote $I_{-1}=[0,1].$
The piecewise polynomial space of degree at most $k$ on grid $\Omega_n$ for $n\ge 0$ is denoted by
\begin{equation}\label{eq:DG-space-Vn}
V_n^k:=\{v: v \in P^k(I_{n}^j),\, \forall \,j=0, \ldots, 2^n-1\}.
\end{equation}
Observing the nested structure
$$V_0^k \subset V_1^k \subset V_2^k \subset V_3^k \subset  \cdots,$$
we can define the multiwavelet subspace $W_n^k$, $n=1, 2, \ldots $ as the orthogonal complement of $V_{n-1}^k$ in $V_{n}^k$ with respect to the $L^2$ inner product on $[0,1]$, i.e.,
\begin{equation*}
V_{n-1}^k \oplus W_n^k=V_{n}^k, \quad W_n^k \perp V_{n-1}^k.
\end{equation*}
By letting $W_0^k:=V_0^k$, we obtain a hierarchical decomposition $V_n^k=\bigoplus_{0 \leq l \leq n} W_l^k$, i.e., MRA of space $V_n^k$.
A set of orthonormal basis can be defined on $W_l^k$ as follows. When $l=0$, the basis $v^0_{i,0}(x)$, $ i=0,\ldots,k$ are the normalized shifted Legendre polynomials in $[0,1]$. When $l>0$, the Alpert's orthonormal multiwavelets \cite{alpert1993class} are employed as the bases and denoted by 
$$v^j_{i,l}(x),\quad i=0,\ldots,k,\quad j=0,\ldots,2^{l-1}-1.$$

We then follow a tensor-product approach to construct the hierarchical finite element space in multi-dimensional space.  Denote $\bl=(l_1,\cdots,l_d)\in\mathbb{N}_0^d$ as the mesh level in a multivariate sense, where $\mathbb{N}_0$  denotes the set of nonnegative integers, we can define the tensor-product mesh grid $\Omega_\bl=\Omega_{l_1}\otimes\cdots\otimes\Omega_{l_d}$ and the corresponding mesh size $h_\bl=(h_{l_1},\cdots,h_{l_d}).$ Based on the grid $\Omega_\bl$, we denote  $I_\bl^\bj=\{\bx:x_m\in(h_mj_m,h_m(j_{m}+1)),m=1,\cdots,d\}$ as an elementary cell, and 
$$\bV_\bl^k:=\{\bv:  \bv \in Q^k(I^{\bj}_{\bl}), \,\,  \bzero \leq \bj  \leq 2^{\bl}-\bone \}= V_{l_1,x_1}^k\times\cdots\times  V_{l_d,x_d}^k$$
as the tensor-product piecewise polynomial space, where $Q^k(I^{\bj}_{\bl})$ represents the collection of polynomials of degree up to $k$ in each dimension on cell $I^{\bj}_{\bl}$. 
If we use equal mesh refinement of size $h_N=2^{-N}$ in each coordinate direction, the  grid and space will be denoted by $\Omega_N$ and $\bV_N^k$, respectively.  
Based on a tensor-product construction, the multi-dimensional increment space can be  defined as
$$\bW_\bl^k=W_{l_1,x_1}^k\times\cdots\times  W_{l_d,x_d}^k.$$
The basis functions in multi-dimensions are defined as
\begin{equation}\label{eq:multidim-basis}
v^\bj_{\bi,\bl}(\bx) := \prod_{m=1}^d v^{j_m}_{i_m,l_m}(x_m),
\end{equation}
for $\bl \in \mathbb{N}_0^d$, $\bj \in B_\bl := \{\bj\in\mathbb{N}_0^d: \,\mathbf{0}\leq\bj\leq\max(2^{\bl-\mathbf{1}}-\mathbf{1},\mathbf{0}) \}$ and $\mathbf{1}\leq\bi\leq \bk+\mathbf{1}$. %The orthonormality of the bases can be easily verified.

Using the notation of $$
|\bl|_1:=\sum_{m=1}^d l_m, \qquad   |\bl|_\infty:=\max_{1\leq m \leq d} l_m.
$$
and  the same component-wise arithmetic operations and relations   as defined in \cite{wang2016elliptic},  we reach the decomposition
\begin{equation}\label{eq:hiere_tp}
\bV_N^k=\bigoplus_{\substack{ |\bl|_\infty \leq N\\\bl \in \mathbb{N}_0^d}} \bW_\bl^k.
\end{equation}
On the other hand, a standard choice of sparse grid   space  \cite{wang2016elliptic, guo2016sparse} is
\begin{equation}
\label{eq:hiere_sg}
\hat{\bV}_N^k=\bigoplus_{\substack{ |\bl|_1 \leq N\\\bl \in \mathbb{N}_0^d}}\bW_\bl^k \subset \bV_N^k.
\end{equation}
We skip the details about the property of the space, but refer the readers to \cite{wang2016elliptic, guo2016sparse}. In Section \ref{sec:dg}, we will describe   the adaptive scheme which adapts a subspace of $\bV_N^k$ according to the numerical solution, hence offering more flexibility and efficiency.

\subsection{Interpolatory multiwavelets}\label{subsec:interp-basis}

Alpert's multiwavelets described in Section \ref{subsec:alpert-basis} are associated with the $L^2$ projection operator. The interpolatory multiwavelets introduced in \cite{tao2019collocation} are constructed based on interpolation operators and also essential for efficient computation of integrals in the DG formulation, especially in high dimensions. In this work, only Lagrange interpolation is considered, while we note that Hermite interpolation can also be used but its implementation is more involved. The details are provided below.

We first define the set of interpolation points on the interval $I=[0,1]$ at zeroth mesh level by $X_0 = \{ x_i \}_{i=0}^M\subset I$. Here, the number of points in $X_0$ is $(M+1)$. We defer the discussion of the relations between $M$ and $k$ to Section \ref{sec:schemeinter}.

The interpolation points at mesh level $n\ge1$, $X_n$ can be obtained correspondingly as
\begin{equation*}
X_n = \{ x_{i,n}^j := 2^{-n}(x_i+j), \quad i=0,\dots,M, \quad j=0,\dots,2^{n}-1 \}.
\end{equation*}
We require the points to be nested, i.e. 
\begin{equation}
\label{nestpts}
X_0 \subset X_1 \subset X_2 \subset X_3 \subset \cdots.
\end{equation}
This can be achieved by requiring  $X_0\subset X_1$.

Given the nodes, we define the basis functions on the zeroth level grid as  Lagrange interpolation polynomials of degree $\le M$ which satisfy the property:
\begin{equation*}
\phi_{i}(x_{i'}) = \delta_{ii'},
\end{equation*}
for $ i,i'=0,\dots,M$. It is easy to see that $\textrm{span} \{ \phi_{i},  i=0,\dots,M \}=V_0^M.  $ 
With the basis function at mesh level zero, we can define the basis functions at mesh level $n\ge1$:
\begin{equation*}
\phi_{i,n}^j := \phi_{i}(2^nx-j), \quad i=0,\dots,M, \quad j=0,\dots,2^n-1, 
\end{equation*}
which form a complete basis set for $V_n^M.$

We now introduce the hierarchical representations and the interpolatory multiwavelets. Define $\tilde{X}_0 := X_0$ and { $\tilde{X}_n := X_n\backslash X_{n-1}$} for $n\ge1$, then we have the decomposition
\begin{equation*}
X_n = \tilde{X}_0 \cup \tilde{X}_1 \cup \cdots \cup \tilde{X}_n .
\end{equation*}
Denote the points in $\tilde{X}_1$ by $\tilde{X}_1=\{ \tilde{x}_i \}_{i=0}^M$. Then the points in $\tilde{X}_n$ for $n\ge1$ can be represented by
\begin{equation*}
\tilde{X}_n = \{ \tilde{x}_{i,n}^j:=2^{-(n-1)}(\tilde{x}_i+j), \quad  i=0,\dots,M, \quad j=0,\dots,2^{n-1}-1 \}.
\end{equation*}

%\begin{equation*}
%    \tilde{W}_1^M = \textrm{span}\{ \psi_{i,l}, \quad i = 0,\dots,P, \quad l=0,\dots,K \}
%\end{equation*}
%where

For notational convenience, we let $\tilde{W}_0^M:=V_0^M.$ The increment function space $\tilde{W}_n^M$ for $n\ge1$ is introduced as a function space    that satisfies
\begin{equation}\label{eq:func-space-sum}
V_n^M = V_{n-1}^M \oplus \tilde{W}_n^M,
\end{equation}
%where $W_0^M:=V_0^M,$ and %Note that $W_n^M$ for $n\ge1$ is not uniquely determined with only the constrain \eqref{eq:func-space-sum}. We pin down a set of basis functions in $W_1^M$ 
and is defined through the multiwavelets
$\psi_{i} \in V_1^M$ that satisfies
\begin{equation*}
\psi_{i}(x_{i'}) = 0, \quad \psi_{i}(\tilde{x}_{i'}) = \delta_{i,i'},
\end{equation*}
for $i,i'=0,\dots,M$. Then  $\tilde{W}_n^M$ is given by
\begin{equation*}
\tilde{W}_n^M = \textrm{span} \{ \psi_{i,n}^j:= \psi_{i}(2^{n-1}x-j), \quad i = 0,\dots,M, \quad j=0,\dots,2^{n-1}-1 \}
\end{equation*}

% The construction above has close connection with interpolation operators.   For a given function $f(x)\in C^{K+1}(I)$, we define $\mathcal{I}^{P,K}_{N}[f]$ as the standard Hermite interpolation on $V^{M}_{N},$ and have the representation
% \begin{align*}% \label{eq:1D_interpolation2}
% 	\mathcal{I}^{P,K}_{N}[f](x) 
% 	% = \sum_{n=0}^{N} \widetilde{\mathcal{I}}^{k}_{n}[f](x)
% 	= \sum_{n=0}^{N} \sum _{j=0}^{\max(2^{n-1}-1,0)}\sum_{l=0}^K \sum_{i=0}^{P} b_{i,l,n}^{j} \psi^{j}_{i,l,n}(x).
% \end{align*}
% Clearly, $(\mathcal{I}^{P,K}_{n}-\mathcal{I}^{P,K}_{n-1})[f](x) \in  \tilde{W}_n^M.$ 
% The algorithm converting between the point values and the derivatives $\{ f^{(l)}(x_{i,n}^j) \}$ to hierarchical coefficients $\{ b_{i,l,n}^{j} \}$ is given in \cite{tao2019collocation}, and by a standard  argument in fast wavelet transform, can be performed in $O(M2^n)$ flops.

The multi-dimensional construction follows similar lines as in Section \ref{subsec:alpert-basis}. We let 
$$\tilde{\bW}_\bl^M=\tilde{W}_{l_1,x_1}^M\times\cdots\times  \tilde{W}_{l_d,x_d}^M,$$
then 
\begin{equation*}
\bV_N^M=\bigoplus_{\substack{ |\bl|_\infty \leq N\\\bl \in \mathbb{N}_0^d}} \tilde{\bW}_\bl^M,
\end{equation*}
while the sparse grid approximation space is
\begin{equation*}
\hat{\bV}_N^M=\bigoplus_{\substack{ |\bl|_1 \leq N\\\bl \in \mathbb{N}_0^d}}\tilde{\bW}_\bl^M.
\end{equation*}
Note that the constructions by Alpert's multiwavelets and the interpolatory multiwavelets deduce the same sparse grid space because of the same nested structure. 
Finally, the interpolation operator in multidimension is defined as $\mathcal{I}^{M}_{N}: C(\Omega)\rightarrow \mathbf{V}^{M}_{N}$:
\begin{align*}% \label{eq:multiD_interpolation}
\mathcal{I}^{M}_{N}[f](\mathbf{x}) 
= \sum_{ \substack{ \abs{\mathbf{n}}_{\infty}\leq N \\ \mathbf{0}\leq \mathbf{j} \leq \max(2^{\mathbf{n}-1}-\mathbf{1},\mathbf{0}) \\ \mathbf{0}\leq \mathbf{i}\leq \mathbf{M}  } } b^{\mathbf{j}}_{\mathbf{i},  \mathbf{n}} \psi^{\mathbf{j}}_{\mathbf{i},  \mathbf{n}} (\mathbf{x}),
\end{align*}
where the multi-dimensional basis functions $\psi^{\mathbf{j}}_{\mathbf{i}, \mathbf{n}} (\mathbf{x})$ are defined in the same approach as \eqref{eq:multidim-basis} by tensor products:
\begin{equation}\label{eq:multidim-basis-interpolation}
\psi^{\mathbf{j}}_{\mathbf{i}, \mathbf{n}} (\mathbf{x}) := \prod_{m=1}^d \psi^{j_m}_{i_m,n_m}(x_m).
\end{equation}
%If the space is switched from $\mathbf{V}^{M}_{N}$ to some subset of $\mathbf{V}^{M}_{N},$ e.g. the sparse grid space $\hat{\bV}_N^M$ or some other subset of $\mathbf{V}^{M}_{N}$ that is dynamically chosen,
For the sparse grid space $\hat{\bV}_N^M$ or any adaptively chosen subspace of $\bV_N^M,$ the interpolation operator, which is denoted by $\mathcal{I}_h^M$ in later sections, can be defined accordingly, by taking only multiwavelet basis functions that belong to that space. The fast algorithms which transform point values at interpolation points to hierarchical coefficients are given in \cite{tao2019collocation}. The detailed formulas of the interpolation points and the associated interpolatory multiwavelets used in this work, we refer readers to \cite{huang2019adaptive,huang2020adaptive}.

\section{Adaptive multiresolution LDG scheme}
\label{sec:dg}
% \todo{need to check if all $h$ has been replaced by $h_N$}

In this section, we present the adaptive multiresolution LDG method for simulating  the HJ equation \eqref{eq:hj}. We start with reviewing the LDG formulation by Yan and Osher in \cite{yan2011local}. Then, by incorporating MRA and multiwavelets introduced in the previous section, we define our scheme.%the sparse grid as well as adaptive multisolution DG methods which are more efficient than the original LDG methods for the HJ equations in high dimensions.    

\subsection{LDG formulation}
\label{sec:semi1}

We consider periodic boundary conditions for simplicity, while the method can be adapted to other non-periodic boundary conditions. 
For illustrative purposes, we first introduce a set of shorthand notation.
Denote by $\Gamma $ the union of the boundaries for all the elements in the partition $\Omega_N$.
The jump  and average of  $q\in L^2(\Gamma)$ are defined as
\begin{align}
[q]=q^- \textbf{n}^- + q^+ \textbf{n}^+, \qquad & \{ q\} = \frac{1}{2} (q^-+q^+),  \nonumber 
%[\textbf{q}] = \textbf{q}^- \cdot \textbf{n}^-  + \textbf{q}^+ \cdot \textbf{n}^+ ,  \qquad & \{ \textbf{q}\} = \frac{1}{2}(\textbf{q}^- + \textbf{q}^+),
\end{align}
where $\textbf{n}$ is the unit normal. `$-$' and `+' represent that the directions of the vector point to interior and exterior at $e$, respectively. Note that $[q]\in [L^2(\Gamma)]^d$, and we  let $[q]_m$ denote the $m$-th component of $[q]$.
% If $e$ is part of the boundary, then we let $[q] = q \textbf{n}$ ($\textbf{n}$ is the outward unit normal) and $\{\textbf{q} \} = \textbf{q}$. 

The key idea in \cite{yan2011local} is to employ the standard LDG methodology, see e.g. \cite{cockburn1998localDG}, to reconstruct the first derivatives of $\phi$, i.e., $\phi_{x_m}$, $m=1,\dots,d$. In particular,  the LDG method  computes two piecewise polynomials $p^1_{m}$ and $p_m^2$, both approximating $\phi_{x_m}$ but using opposite one-sided numerical fluxes; that is, given $\phi_h$ we seek $p_m^\tau$, $m=1,\ldots,d$, $\tau=1,\,2$ in $\bV$ such that for all $w_h\in\bV$
\begin{align}\label{eq:ldg_deri}
\int_{\Omega}p_m^\tau w_h\,d\bx = -\int_{\Omega} \phi_{h} (w_h)_{x_m} d \mathbf{x}+\sum_{e \in \Gamma} \int_{e} \widehat{\phi}^\tau_h\left[w_{h}\right]_m d s
\end{align}
where the numerical fluxes are defined as
$$
\widehat{\phi}^1_h =  \{\phi_h\} + \frac12[\phi_h]_m,\quad
\widehat{\phi}^2_h =  \{\phi_h\} - \frac12[\phi_h]_m.
$$
Note that $p_m^1$ and $p_m^2$ carry the information of $\phi_{x_m}$ from opposite directions. Hence, when the solution is smooth,  $p_m^1$ and $p_m^2$ are almost identical, while if the solution involves nonsmooth corners, then   $p_m^1$ and $p_m^2$ can be very different.

%and $\mathbf{e}_m$ is the $m$-th standard unit vector. 
Then the semi-discrete scheme for solving \eqref{eq:hj} is defined as follows: seek $\phi_h\in \bV$ such that, for all $v\in\bV$,
\begin{equation}
\label{eq:ldg_hj}
\int_{\Omega}(\phi_h)_t v\ d\bx +\int_{\Omega} \widehat{H}(p^1_1,p^2_1,p^1_2,p^2_2,\ldots,p^1_d,p^2_d)v\,d\bx = 0,
\end{equation}
where $\widehat{H}$ denotes a monotone numerical Hamiltonian that approximates $H$, and $p_m^\tau$, $\tau=1,2$, $m=1,\ldots,d$ are given in \eqref{eq:ldg_deri}. In the simulations, we employ the following global Lax-Friedrichs Hamiltonian
$$
\widehat{H}(p^1_1,p^2_1,p^1_2,p^2_2,\ldots,p^1_d,p^2_d) = H(\bar{p}_1,\bar{p}_2,\ldots,\bar{p}_d) - \sum_{m=1}^d \frac{\alpha_m}{2}\left(p^2_m-p^1_m\right),
$$
where $\bar{p}_m = \frac12(p^1_m + p^2_m)$ and 
$$\alpha_m=\max_{q_1,\ldots,q_d} \left|\frac{\partial H(q_1,\ldots,q_d)}{\partial q_m}\right|$$ 
with the maximum being taken over the whole domain.

Depending on the choice of space $\bV$, we obtain several LDG methods for \eqref{eq:hj} with distinct properties.  If   $\bV=\bV^k_N,$ we recover the full grid LDG scheme in \cite{yan2011local} on tensor-product meshes. If $\bV=\hat{\bV}^k_N,$ then we obtain  the sparse grid LDG method. If $\bV$ is chosen adaptively, we have the adaptive sparse grid scheme. Noteworthy,  besides the LDG formulation,  we can employ other  DG formulations as well, such as the direct DG method \cite{cheng2014hj} and the indirect DG methods \cite{hu1999discontinuous,guo2011local}. The LDG formulation used is comparatively simpler to implement under the MRA framework. 

If the Hamiltonian $H$ is linear, then the HJ equation \eqref{eq:hj} degenerates to a transport equation with constant coefficients, and the formulation \eqref{eq:ldg_hj} together with \eqref{eq:ldg_deri} is nothing but a standard upwind DG scheme. The results established in \cite{guo2016sparse} can be adapted directly to the linear HJ equation that if the solutions is adequately smooth in terms of the mixed norm, then the sparse grid DG method using space $\hat{\bV}_n^k$ is convergent of order $k+1$ with a polylogarithmic factor.

\subsection{Semi-discrete scheme with multiresolution interpolation}
\label{sec:schemeinter}

%To evaluate the integrals over elements and edges more efficiently with a cost proportional to the DoF of the underlying finite element space, 

 For nonlinear problems, one major difficulty of implementation of formulation \eqref{eq:ldg_hj} is to compute the volume integral efficiently and accurately, especially in high dimensions. Naive implementation of numerical quadratures is inefficient due to the hierarchical structure of  multiwavelets. To address the challenge, we follow the idea in \cite{shen2010efficient,huang2019adaptive} and interpolate the numerical Hamiltonian $\widehat{H}$ by using the multiresolution Lagrange interpolation  discussed in Section \ref{subsec:interp-basis}. In particular, we have the following modified formulation with interpolation. We find $\phi_h\in\bV$  so that for all $v\in\bV$ 
 \begin{equation}
 \label{eq:ldg_hj_inter}
 \int_{\Omega}(\phi_h)_t v\ d\bx +\int_{\Omega} \mathcal{I}_h^M\left(\widehat{H}(p^1_1,p^2_1,p^1_2,p^2_2,\ldots,p^1_d,p^2_d)\right)v\,d\bx = 0.
 \end{equation}
By doing so, not only can we apply the unidirectional principle to facilitate the computation, but also fast algorithms can be utilized to improve efficiency. We omit the details regarding the fast algorithm and refer readers to \cite{huang2019adaptive}.

%Note that the interpolation procedure in the scheme formulation plays a role as a high order sparse grid numerical quadrature. It is also known that a numerical quadrature with adequate accuracy has to be used for DG methods to ensure stability when solving nonlinear conservation laws \cite{}. We have similar observation here for solving nonlinear HJ equations, i.e., the interpolation must have sufficient accuracy for stability. 
Note that the interpolation procedure in the scheme formulation plays a role as a high order MRA numerical quadrature. There exist two types of points, namely inner  points and interface points \cite{tao2019collocation}. It is observed that the schemes using the interface points are more stable than those using the inner points with the same order accuracy (see, e.g. \cite{huang2019adaptive}), and hence we choose to use the interface points in the simulations. We also remark if the Hamiltonian $H$ is not smooth and $k>1$, then to ensure stability, one must employ a very high order quadrature, i.e. large $M$, for accurate computation of the volume integral. This is ascribed to the fact that the large quadrature error due to nonsmoothness of $H$ may pollute the numerical viscosity and lead to instability. This drawback is observed in \cite{yan2011local}, and the authors further coupled a nonlinear limiter to restore stability. In this paper, we propose to properly regularize  the Hamiltonian so that the interpolation $\mathcal{I}_h^M$ is adequately accurate for stability. The details will be presented in Section \ref{sec:numer}. Another possible approach is to add artificial viscosity, as done for solving conservation laws \cite{huang2019adaptive}.

%In the one-dimensional setting of formulation \eqref{eq:ldg_hj_inter}, they indeed correspond to the famous open and closed Newton-Cotes quadrature formula, respectively.

To preserve the accuracy of the original DG scheme, the interpolation operator $\mathcal{I}^M_h(\cdot)$ needs to reach certain accuracy. Following \cite{DG4}, we can write the DG scheme with interpolation \eqref{eq:ldg_hj_inter} into the semi-discrete form as
\begin{equation}
\label{eq:ode}
	\frac{d\phi_h}{dt} = L_h(\phi_h).
\end{equation}
Here $L_h(\cdot)$ is an operator onto $\bV$ and is a discrete approximation of $-H(\nabla\phi)$ which satisfies
\begin{equation}
 \int_{\Omega}L_h(\phi_h) v\ d\bx +\int_{\Omega} \mathcal{I}_h^M\left(\widehat{H}(p^1_1,p^2_1,p^1_2,p^2_2,\ldots,p^1_d,p^2_d)\right)v\,d\bx = 0.
\end{equation}
for all $v\in\bV$ with $p_m^\tau$, $\tau=1,2$, $m=1,\ldots,d$ determined by \eqref{eq:ldg_deri}. 
Using  similar   techniques as in \cite{DG4,huang2017quadrature}, we have the following proposition on local truncation error of the sparse grid method with $\bV=\hat{\bV}^k_N.$ The proof is omitted for brevity.
\begin{prop}[Local truncation error analysis]%[Accuracy of semi-discrete sparse grid DG scheme with interpolation]
	\label{prop:accurate-interp}
	%Assume that the sparse DG finite element space has polynomials up to degree $k,$ 
	If the interpolation operator $\mathcal{I}_h^M$ in \eqref{eq:ldg_hj_inter} has the accuracy of order $\abs{\log_2h_N}^{d}h_N^{k+1}$ for sufficiently smooth functions, then the local truncation error of the semi-discrete DG scheme with interpolation \eqref{eq:ldg_hj_inter} is of order $\abs{\log_2h_N}^{d}h_N^{k+1}$. To be more precise, for sufficiently smooth Hamiltonian $H$ and function $\phi$, the sparse grid DG method with interpolation \eqref{eq:ldg_hj_inter} has the truncation error:
	\begin{equation}\label{eq:truncation-sparse}
	\norm{L_h(\phi) + H(\nabla\phi) )}_{L^2(\Omega)}\le C \abs{\log_2h_N}^{d}h_N^{k+1}.
	\end{equation}
	Here, we use  $C$ to denote a generic constant that may depend on the solution $u$, but does not depend on $N.$
\end{prop}

The proposition indicates that, to  preserve the order accuracy of the original scheme, we should use $M \ge k.$ Hence, in the simulation we let $M \ge k$. Meanwhile, many Hamiltonians are non-smooth functions, and indeed we  observe numerically that we need  $M>k$ for those cases. This will be further discussed in Section \ref{sec:numer}.
%	
%	For example, if we take piecewise linear polynomials for the DG space, then it is required to apply linear interpolation operator to treat the nonlinear terms. \todo{From our numerical tests, however, this seems that it is not a necessary condition. To reach the desired convergence rate, one only needs to take $M\ge k$.}

% In Proposition \ref{prop:accurate-interp}, we only estimate the local truncation error, and this is far from a rigorous error estimate that takes into account  stability. Unlike the scheme with the symmetric bilinear form $B(u_h,v)$ as in Theorem \ref{thm:stable}, the symmetry is lost in the interpolated bilinear form $\tilde{B}(u_h, v).$  Hence, energy stability is not automatic. In numerical experiments, we observe that the sparse grid DG method with Lagrange interpolation with only inner interface points is unstable for polynomials of high degrees (see the numerical results in Table \ref{tab:smooth-sparse-2D-inner-pt} in Section 4). With the interpolation points at the {\emph{interface}}, the sparse grid DG scheme is   stable and yields satisfactory convergence rate (see Table \ref{tab:smooth-sparse-2D-interface-pt} in Section 4).

%\subsection{Time stepping and adaptivity}
%\label{subsec:timeadapt}

For time discretization, we employ the third order strong-stability-preserving Runge-Kutta (RK) scheme \cite{shu1988jcp} to advance the semi-discrete scheme \eqref{eq:ode}. 
%The reason why we prefer the one-step RK method to a multistep method is that the maximum allowed time step size from the CFL restriction may change with the adaptive mesh in different time steps. This would result in additional computational cost in extrapolation or interpolation between different time steps for the multistep methods. %If one sticks to apply the fixed time step determined by the finest grid in the multistep method, it will be computationally expensive for coarse mesh.
The adaptive procedure follows the technique developed in \cite{bokanowski2013adaptive,guo2017adaptive} to determine the space $\bV$ that dynamically evolves over time. The details are omitted for brevity. The main idea is that in light of the distinguished property of multiwavelets,  we keep track of multiwavelet coefficients, i.e. $L^2$ norms of $\phi_h$, as an error indicator for refining and coarsening, aiming to efficiently capture the viscosity solution of \eqref{eq:hj} which may develop discontinuous derivatives. 
\section{Numerical examples}
\label{sec:numer}
In this section, we present a collection of numerical examples to demonstrate the performance of the proposed adaptive sparse grid LDG method for solving the HJ equation. We consider numerical examples  up to $d=4$ with smooth and nonsmooth Hamiltonian, and with smooth and nonsmooth viscosity solutions.
Noteworthy, we may need to tune $M$ for optimal performance. In particular, we observe that for some numerical tests, we can simply take $M=k$ to achieve satisfactory results and maintain the original accuracy of the DG method, while for the some other tests, we may need to take larger $M$ to ensure good performance. In all numerical simulations, the value of $M$ is taken between $k$ and $k+2.$

%%%%%%%%%%%%%%%%%%%%%%%%%%%%%%%%%%%%%%%%%%%%%%

\begin{exa}\label{ex:burgers}
Consider the following Burgers' equation in $d$-dimension
\begin{equation}
\label{eq:burgers}
\begin{cases}\displaystyle
\phi_t + \frac12\left(\sum_{m=1}^{d}\phi_{x_m}\right)^2 = 0,\quad \bx\in[0,1]^d,\\
\displaystyle\phi(\bx,0) = -\frac{1}{2\pi}\cos\left(2\pi\sum_{m=1}^d x_m\right),
\end{cases}
\end{equation}
with periodic boundary conditions.
\end{exa}

At $T=0.01$ for $d=2$ and $T=0.005$ for $d=3$, the solutions are still smooth, and we summarize the convergence study of the sparse grid method with $\hat{\bV}_N^k$ in Tables \ref{tb:burgers_sparse}-\ref{tb:burgers_sparse3d}, including the L$^2$ errors and the associated order of accuracy, with various configurations of $k$ and $M$. It is observed that larger $M$  leads to smaller error magnitude as expected. Slight order reduction is observed for $d=2$, and it becomes more severe for $d=3$. This is because, as time evolves, the viscosity solution of \eqref{eq:burgers} develops larger and larger mixed derivatives especially in high dimensions. Hence, it may not be optimal to use the sparse grid space $\hat{\bV}_N^k$ for approximating the viscosity solution. In Tables \ref{tb:burgers_adap_d2}-\ref{tb:burgers_adap_d3} we report the convergence study for the adaptive method for $d=2,3$, respectively. In particular, by fixing the maximum mesh level $N=7$, two rates of convergence are calculated \cite{bokanowski2013adaptive}. The first one is with  respect  to the error threshold:
$$R_{\epsilon_{l}}=\frac{\log \left(e_{l-1} / e_{l}\right)}{\log \left(\epsilon_{l-1} / \epsilon_{l}\right)},$$
and the second is with respect to DoF:
$$R_{\mathrm{DoF}_{l}}=\frac{\log \left(e_{l-1} / e_{l}\right)}{\log \left(\mathrm{DoF}_{l} / \mathrm{DoF}_{l-1}\right)}.$$
For the full grid counterpart, we have $R_{\mathrm{DoF}} = (k+1)/d$ for smooth solutions. It is observed that $R_{\epsilon}<1$, which is similar to the Burgers' equation \cite{huang2019adaptive}. Furthermore, using larger $k$ is beneficial, as the method with larger $k$ requires less  DoF to attain a certain level of accuracy.

\begin{figure}[h!]
	\centering
	\subfigure[]{\includegraphics[height=50mm]{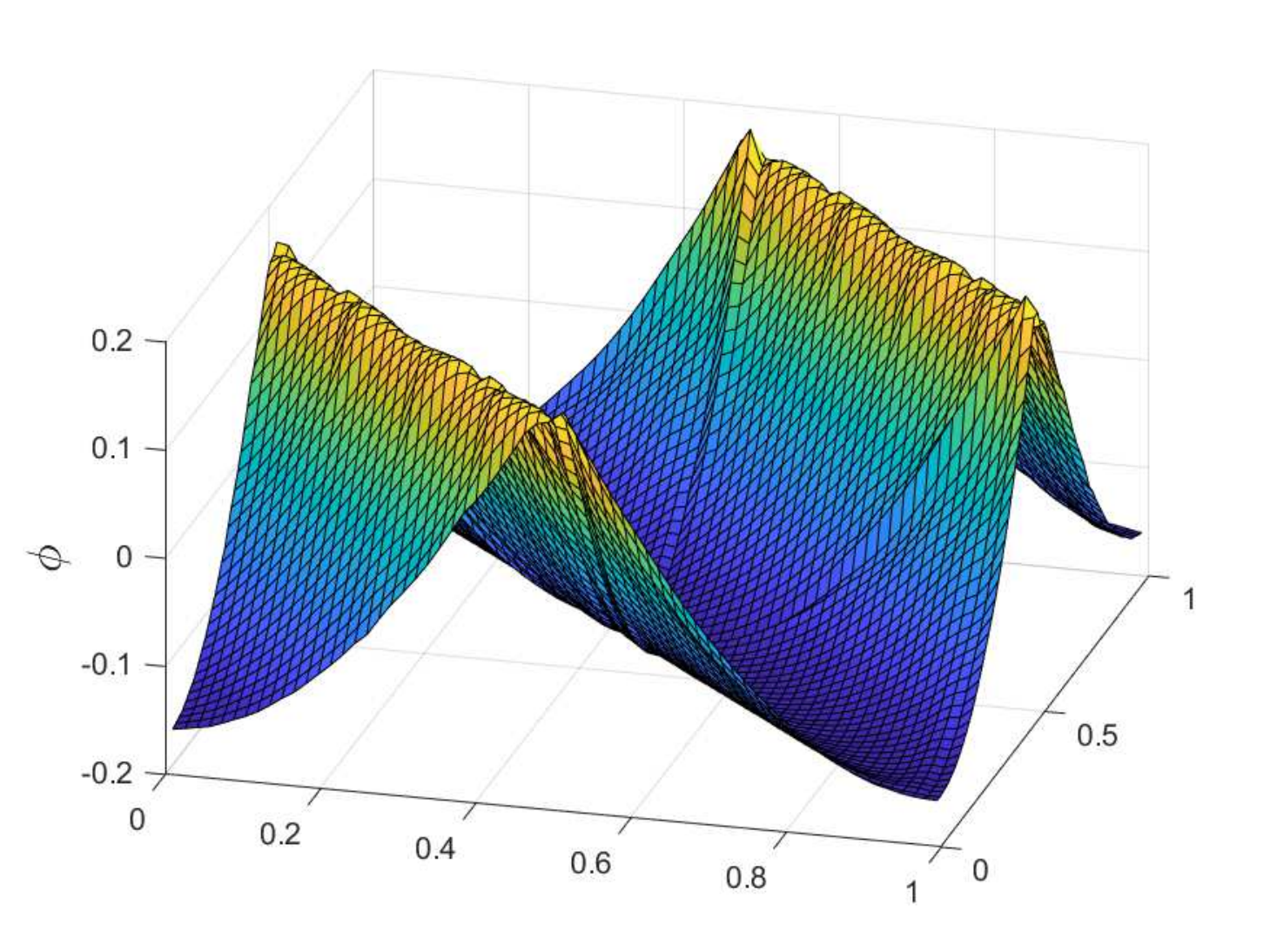}}
	\subfigure[]{\includegraphics[height=50mm]{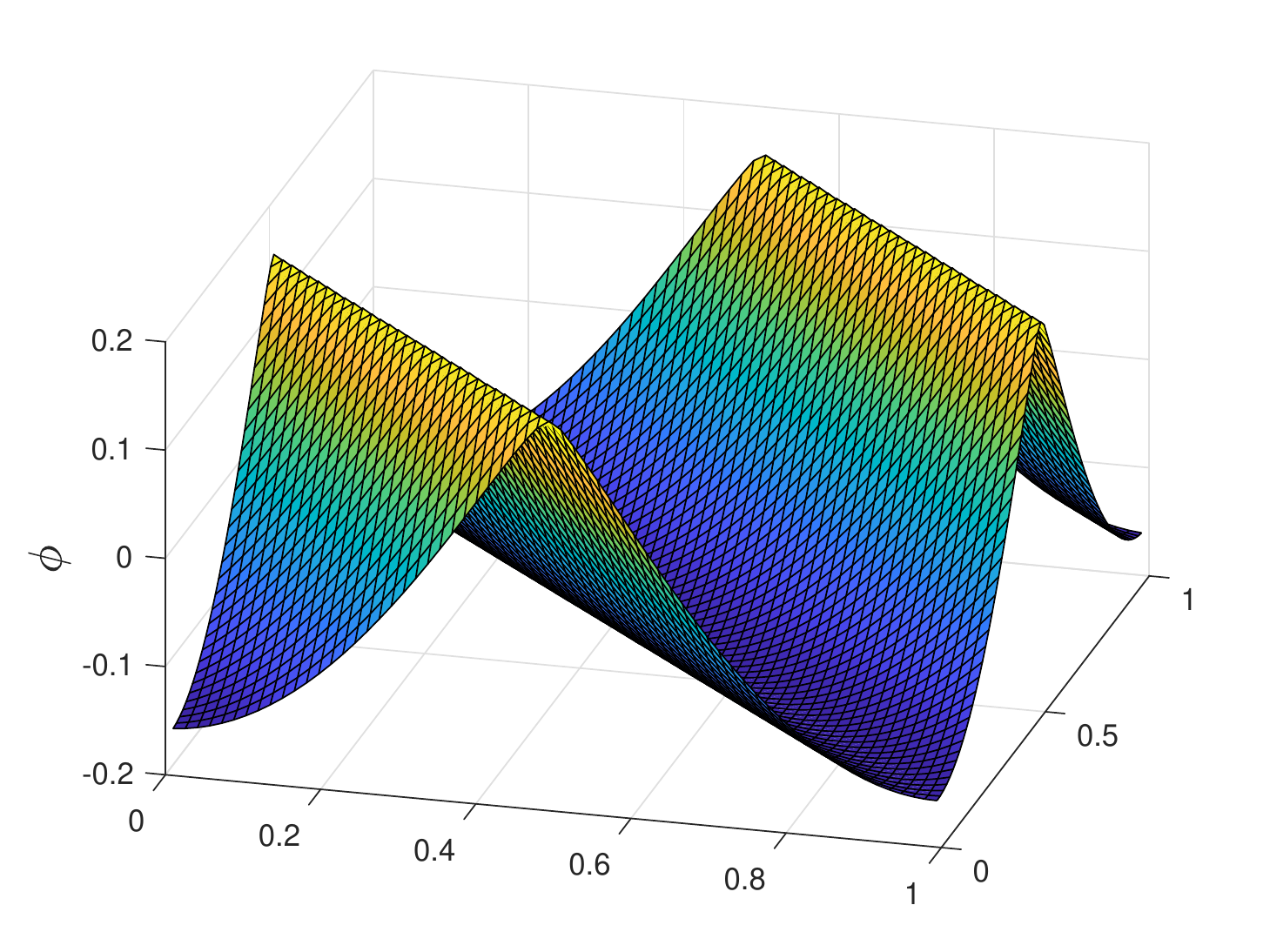}}
	\subfigure[]{\includegraphics[height=50mm]{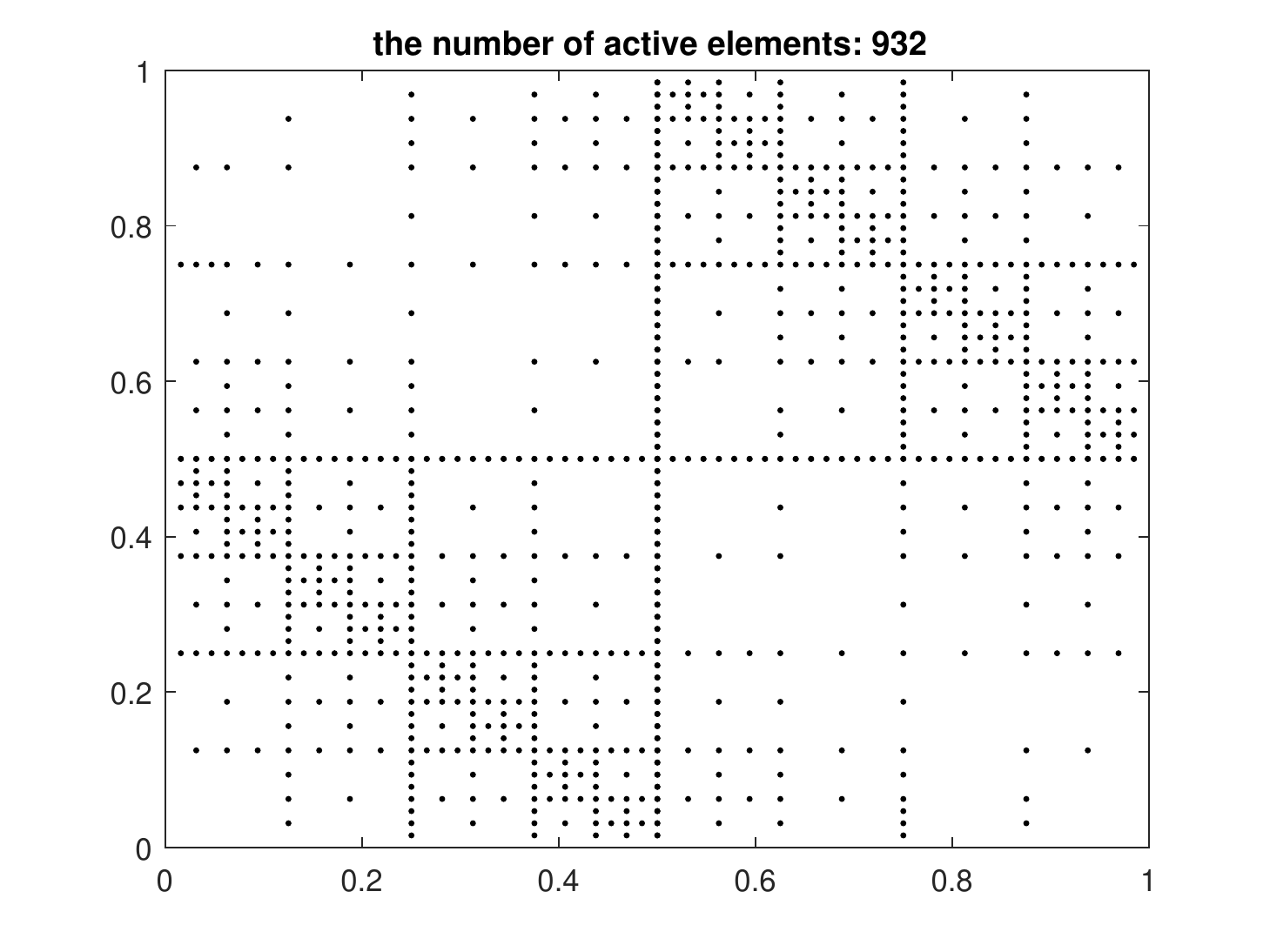}}
	\caption{Example \ref{ex:burgers}, $d=2$. $k=2$, $M=2$. $T=0.04$. $N=6$. $\varepsilon$=$10^{-5}$. (a) Numerical solution by sparse grids. (b) Numerical solutions by adaptive sparse grid. (c) Active elements.  \label{fig:burgers}}
\end{figure}

\begin{figure}[h!]
	\centering
	\subfigure[]{\includegraphics[height=50mm]{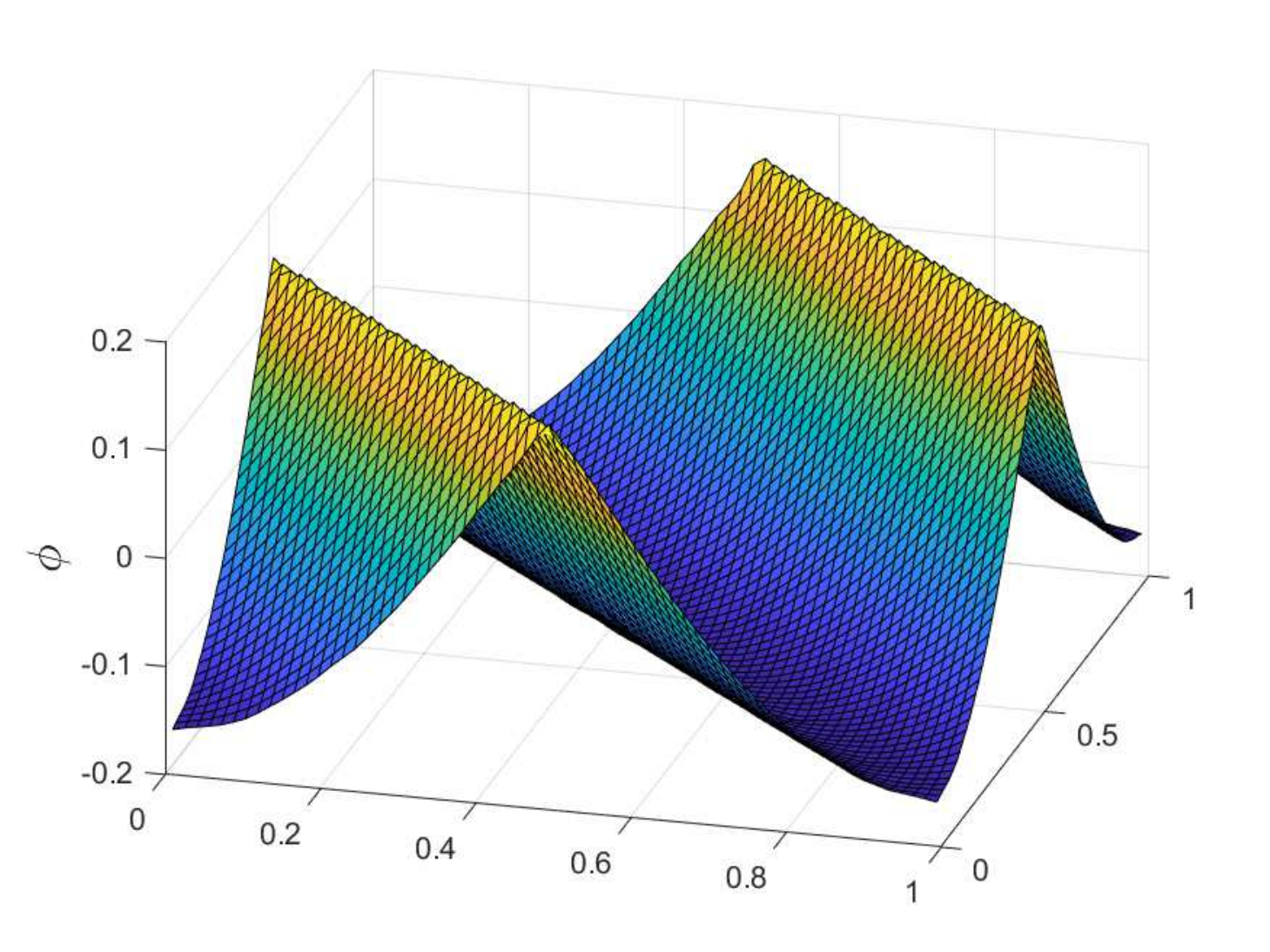}}
	\subfigure[]{\includegraphics[height=50mm]{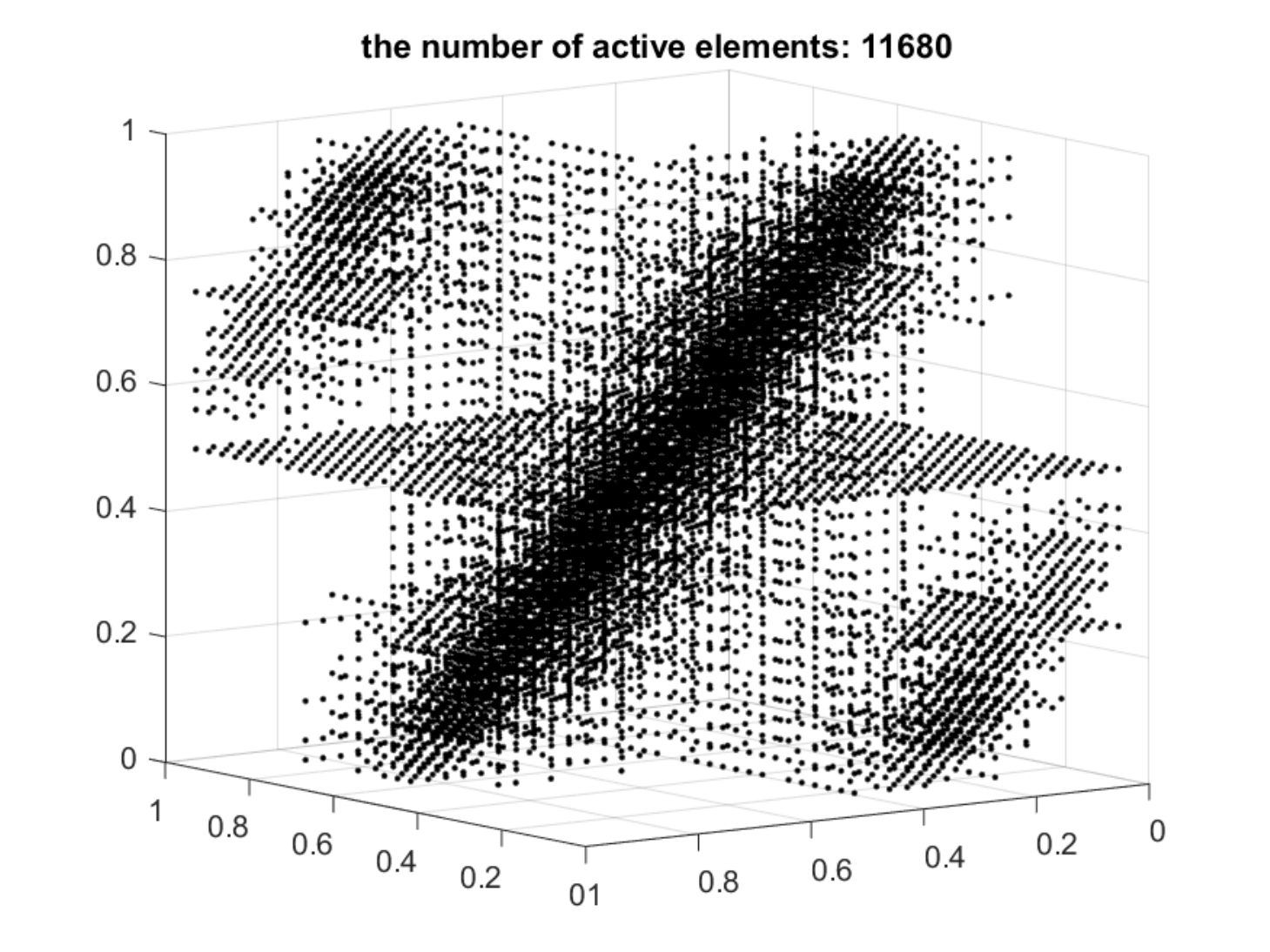}}
	\caption{Example \ref{ex:burgers}, $d=3$. $k=2$, $M=3$. $T=0.02$. $N=6$. $\varepsilon$=$10^{-5}$. (a) 2D-cuts of the numerical solution at $x_3=0$. (b) Active elements.  \label{fig:burgers_3d}}
\end{figure}

At $T=0.04$ for $d=2$ and $T=0.02$ for $d=3$, the viscosity solutions have developed discontinuous derivatives. In Figure \ref{fig:burgers}, we plot the solution profiles computed by the sparse grid method and the adaptive method for $d=2$. We set $k=2$, $M=2$ and $\epsilon = 10^{-5}$, and the maximum mesh level $N=6$ for the adaptive method. It is observed the sparse grid method is able to capture the main structure of the solution, but severe oscillations appear due to lack of mesh resolution around the corners, while the adaptive method is able to effectively capture the viscosity solution by adding more DoF in the nonsmooth region. Hence, the sparse grid method with fixed space $\hat{\bV}_N^k$ cannot reliably approximate the nonsmooth viscosity solution. Afterwards, we will only focus on the performance of adaptive method. For $d=3$, we set $k=2$, $M=3$ and $\epsilon = 10^{-5}$, and the maximum mesh level $N=6$, and plot the results generated by the adaptive method in Figure \ref{fig:burgers_3d}, including the 2D cuts of the solution and the associated active elements at final time. Similar results to $d=2$ are observed.

% For the case of $k=3$ and $M=3$, the sparse version is nonlinearly unstable (it is actually linearly stable. For linear problem, the interpolation recovery the original polynomial and hence the method is same as the original sparse grid DG method, which is stable.) Furthermore, the full grid as well as the adaptive version are observed to be stable. This observation is consistent to the case of simulating hyperbolic conservation laws. 

%
%\begin{table}[!hbp]
%	\centering
%	\caption{Example \ref{ex:burgers}, $d=2$. Sparse grid. $T=0.01$.}
%	\label{tb:burgers_sparse}
%	\begin{tabular}{|c|c|c|c|c|c|c|}
%		\hline
%		 $N$ & L$^1$-error & order & L$^2$-error & order & L$^{\infty}$-error & order \\
%		\hline
%			&	\multicolumn{6}{c|}{$k=1$}\\\hline
%3	&	2.06E-2	&	0.00	&	2.63E-2	&	0.00	&	9.25E-2	&	0.00	\\
%4	&	5.99E-3	&	1.78	&	7.87E-3	&	1.74	&	2.71E-2	&	1.77	\\
%5	&	2.87E-3	&	1.06	&	3.82E-3	&	1.04	&	1.36E-2	&	1.00	\\
%6	&	1.35E-3	&	1.09	&	1.77E-3	&	1.11	&	8.21E-3	&	0.72	\\
%7	&	6.21E-4	&	1.12	&	7.95E-4	&	1.16	&	2.91E-3	&	1.50	\\\hline
%			&	\multicolumn{6}{c|}{$k=2$}\\\hline
%3&	4.52E-3&	-	&5.68E-3	&-	&1.55E-2&	- \\ 
%4&	8.36E-4&	2.43&	1.21E-3&	2.23&	4.32E-3&	1.85\\
%5&	2.24E-4&	1.90&	3.10E-4&	1.97&	1.30E-3&	1.73\\
%6&	3.64E-5&	2.62&	5.28E-5&	2.56&	3.46E-4&	1.91\\
%7&	6.41E-6&	2.51&	9.58E-6&	2.46&	5.97E-5&	2.53\\
%\hline
%
%
%	\end{tabular}
%\end{table}

\begin{table}[htbp]
	\centering
\caption{Example \ref{ex:burgers}, $d=2$. Sparse grid. $T=0.01$.}
	\label{tb:burgers_sparse}%
	\begin{tabular}{c|c|c|c|c|c|c|c}
		\hline
		\multirow{2}[4]{*}{} & \multirow{2}[4]{*}{$N$} & \multicolumn{2}{c|}{$M=1$} & \multicolumn{2}{c|}{$M=2$} & \multicolumn{2}{c}{$M=3$} \\
		\cline{3-8}          &       & $L_2$ error & order & $L_2$ error & order & $L_2$ error & order \\
		\hline
		\multirow{5}[3]{3em}{$k=1$} 
&3	&	2.63E-02	&	--	&		1.99E-02	&	--	&	1.99E-02	&	--	\\
&4	&	7.87E-03	&	1.74	&		5.88E-03	&	1.76	&	5.75E-03	&	1.79	\\
&5	&	3.82E-03	&	1.04	&		2.42E-03	&	1.28	&	2.24E-03	&	1.36	\\
&6	&	1.77E-03	&	1.11	&		8.27E-04	&	1.55	&	6.64E-04	&	1.76	\\
&7	&	7.95E-04	&	1.16	&		3.47E-04	&	1.25	&	2.10E-04	&	1.66	\\
		\hline
		\hline
		\multirow{2}[4]{*}{} & \multirow{2}[4]{*}{$N$} & \multicolumn{2}{c|}{$M=2$} & \multicolumn{2}{c|}{$M=3$} & \multicolumn{2}{c}{$M=4$} \\
		\cline{3-8}          &       &  $L_2$ error & order &  $L_2$ error & order &  $L_2$ error & order \\
		\hline
		\multirow{5}[3]{3em}{$k=2$} 
&3	&	5.68E-03	&	--	&	2.84E-03	& --	&	2.81E-03	&	--	\\
&4	&	1.21E-03	&	2.23	&	3.75E-04	&	2.92	&	3.66E-04	&	2.94	\\
&5	&	3.10E-04	&	1.97	&	1.42E-04	&	1.40	&	1.41E-04	&	1.37	\\
&6	&	5.28E-05	&	2.56	&	1.97E-05	&	2.84	&	1.96E-05	&	2.85	\\
&7	&	9.58E-06	&	2.46	&	4.36E-06	&	2.18	&	4.32E-06	&	2.18	\\

		\hline  
		\hline  
		\multirow{2}[4]{*}{} & \multirow{2}[4]{*}{$N$} & \multicolumn{2}{c|}{$M=3$} & \multicolumn{2}{c|}{$M=4$} & \multicolumn{2}{c}{$M=5$} \\
		\cline{3-8}          &       &  $L_2$ error & order &  $L_2$ error & order &  $L_2$ error & order \\
		\hline
		\multirow{5}[3]{3em}{$k=3$} 
&3	&	1.14E-03	&	--	&	6.28E-04	&--	&	6.10E-04	&	--	\\
&4	&	1.68E-04	&	2.76	&	6.86E-05	&	3.19	&	6.58E-05	&	3.21	\\
&5	&	2.59E-05	&	2.70	&	1.45E-05	&	2.24	&	1.44E-05	&	2.20	\\
&6	&	2.10E-06	&	3.63	&	7.94E-07	&	4.19	&	7.84E-07	&	4.20	\\
&7	&	2.77E-07	&	2.92	&	1.50E-07	&	2.41	&	1.49E-07	&	2.39	\\

		\hline    
	\end{tabular}%  
\end{table}%

\begin{table}[htbp]
	\centering
	\caption{Example \ref{ex:burgers}, $d=3$. Sparse grid. $T=0.005$.}
	\label{tb:burgers_sparse3d}%
	\begin{tabular}{c|c|c|c|c|c|c|c}
		\hline
		\multirow{2}[4]{*}{} & \multirow{2}[4]{*}{$N$} & \multicolumn{2}{c|}{$M=1$} & \multicolumn{2}{c|}{$M=2$} & \multicolumn{2}{c}{$M=3$} \\
		\cline{3-8}          &       & $L_2$ error & order & $L_2$ error & order & $L_2$ error & order \\
		\hline
		\multirow{6}[3]{3em}{$k=1$} 
&	3	&	3.34E-02	&	--	&	9.22E-03	&	--	&	7.46E-03	&	--	\\
&	4	&	1.71E-02	&	0.96	&	3.69E-03	&	1.32	&	3.24E-03	&	1.20	\\
&	5	&	6.93E-03	&	1.31	&	1.32E-03	&	1.49	&	1.21E-03	&	1.43	\\
&	6	&	2.07E-03	&	1.75	&	5.49E-04	&	1.26	&	5.28E-04	&	1.19	\\
&	7	&	8.07E-04	&	1.36	&	1.57E-04	&	1.80	&	1.54E-04	&	1.78	\\
&	8	&	1.91E-04	&	2.08	&	4.09E-05	&	1.94	&	3.98E-05	&	1.95	\\

		\hline
		\hline
		\multirow{2}[4]{*}{} & \multirow{2}[4]{*}{$N$} & \multicolumn{2}{c|}{$M=2$} & \multicolumn{2}{c|}{$M=3$} & \multicolumn{2}{c}{$M=4$} \\
		\cline{3-8}          &       &  $L_2$ error & order &  $L_2$ error & order &  $L_2$ error & order \\
		\hline
		\multirow{6}[3]{3em}{$k=2$} 
&	3	&	3.34E-02	&	--	&	9.22E-03	&	--	&	7.46E-03	&	--	\\
&	4	&	1.71E-02	&	0.96	&	3.69E-03	&	1.32	&	3.24E-04	&	1.20	\\
&	5	&	6.93E-03	&	1.31	&	1.32E-03	&	1.49	&	1.21E-03	&	1.43	\\
&	6	&	2.07E-03	&	1.75	&	5.49E-04	&	1.26	&	5.28E-04	&	1.19	\\
&	7	&	8.07E-04	&	1.36	&	1.57E-04	&	1.80	&	1.54E-04	&	1.78	\\
&	8	&	1.91E-04	&	2.08	&	4.09E-05	&	1.94	&	3.98E-05	&	1.95	\\

		\hline  
		\hline  
		\multirow{2}[4]{*}{} & \multirow{2}[4]{*}{$N$} & \multicolumn{2}{c|}{$M=3$} & \multicolumn{2}{c|}{$M=4$} & \multicolumn{2}{c}{$M=5$} \\
		\cline{3-8}          &       &  $L_2$ error & order &  $L_2$ error & order &  $L_2$ error & order \\
		\hline
		\multirow{6}[3]{3em}{$k=3$} 
&	3	&	8.86E-03	&	--	&	3.56E-03	&	--	&	2.48E-03	&	--	\\
&	4	&	2.97E-03	&	1.58	&	1.10E-03	&	1.70	&	8.62E-04	&	1.53	\\
&	5	&	9.97E-04	&	1.57	&	3.64E-04	&	1.59	&	2.93E-04	&	1.56	\\
&	6	&	3.08E-04	&	1.70	&	9.78E-05	&	1.90	&	8.57E-05	&	1.77	\\
&	7	&	6.49E-05	&	2.24	&	2.17E-05	&	2.17	&	2.04E-05	&	2.07	\\
&	8	&	1.44E-05	&	2.17	&	5.02E-06	&	2.12	&	4.85E-06	&	2.07	\\
\hline

	\end{tabular}%  
\end{table}%

\begin{table}[!hbp]
	\centering
	\caption{Example \ref{ex:burgers}, $d=2$. Adaptive sparse grid. $T=0.01$. $M=k$.}
	\label{tb:burgers_adap_d2}
	\begin{tabular}{c|c|c|c|c|c}
		\hline
		& $\epsilon$ & DoF & L$^2$-error  & $R_{\epsilon}$ & $R_{\textrm{DoF}}$\\
		\hline
		\multirow{4}{3em}{$k=1$}
&	1.00E-03	&	448	&	1.56E-03	&	--	&	--	\\
&	1.00E-04	&	1376	&	6.92E-04	&	0.35	&	0.73	\\
&	1.00E-05	&	3520	&	2.55E-04	&	0.43	&	1.06	\\
&	1.00E-06	&	10240	&	2.26E-05	&	1.05	&	2.27	\\
&	1.00E-07	&	18688	&	1.05E-05	&	0.33	&	1.28	\\
 
		\hline\hline
		\multirow{4}{3em}{$k=2$}
&	1.00E-03	&	270	&	6.67E-04	&	--	&	--	\\
&	1.00E-04	&	720	&	3.10E-04	&	0.33	&	0.78	\\
&	1.00E-05	&	1548	&	4.93E-05	&	0.80	&	2.40	\\
&	1.00E-06	&	3492	&	1.34E-05	&	0.57	&	1.60	\\
&	1.00E-07	&	7704	&	1.83E-06	&	0.86	&	2.51	\\

		\hline\hline
		\multirow{4}{3em}{$k=3$}
&	1.00E-03	&	192	&	1.14E-03	&	--	&	--	\\
&	1.00E-04	&	480	&	1.04E-04	&	1.04	&	2.62	\\
&	1.00E-05	&	896	&	3.14E-05	&	0.52	&	1.92	\\
&	1.00E-06	&	1856	&	7.14E-06	&	0.64	&	2.03	\\
&	1.00E-07	&	3136	&	7.07E-07	&	1.00	&	4.41	\\
		\hline
	\end{tabular}
\end{table}

\begin{table}[!hbp]
	\centering
	\caption{Example \ref{ex:burgers}, $d=3$. Adaptive sparse grid. $T=0.005$. $M=k$.}
	\label{tb:burgers_adap_d3}
	\begin{tabular}{c|c|c|c|c|c}
		\hline
		& $\epsilon$ & DoF & L$^2$-error  & $R_{\epsilon}$ & $R_{\textrm{DoF}}$\\
		\hline
		\multirow{5}{3em}{$k=1$}
&1.00E-03	&	2432	&	7.87E-03	&	--	&	--	\\
&1.00E-04	&	14864	&	3.03E-03	&	0.41	&	0.53	\\
&1.00E-05	&	44656	&	1.17E-03	&	0.41	&	0.87	\\
&1.00E-06	&	152176	&	3.25E-04	&	0.56	&	1.04		\\
&1.00E-07	&	380976	&	9.05E-05	&	0.56	&	1.39		\\

		\hline\hline
		\multirow{5}{3em}{$k=2$}
&1.00E-03	&	2646	&	5.84E-03	&	--	&	--	\\
&1.00E-04	&	8208	&	9.84E-04	&	0.77	&	1.57	\\
&1.00E-05	&	21816	&	1.96E-04	&	0.70	&	1.65	\\
&1.00E-06	&	55404	&	6.11E-05	&	0.51	&	1.25	\\
&1.00E-07	&	133569	&	1.33E-05	&	0.66	&	1.74	\\
		
		\hline\hline
		\multirow{5}{3em}{$k=3$}
&1.00E-03	&	2048	&	2.01E-03	&	--	&	--	\\
&1.00E-04	&	6400	&	5.01E-04	&	0.60	&	1.22	\\
&1.00E-05	&	16384	&	9.26E-05	&	0.73	&	1.80	\\
&1.00E-06	&	35584	&	2.30E-05	&	0.60	&	1.79	\\
&1.00E-07	&	99840	&	2.82E-06	&	0.91	&	2.04	\\
		\hline
	\end{tabular}
\end{table}

\begin{exa}\label{ex:cos}
	Consider the following HJ equation with a nonconvex Hamiltonian
	\begin{equation}
	\begin{cases}\displaystyle
	\phi_t - \cos\left(\sum_{m=1}^{d}\phi_{x_m} + 1\right) = 0,\quad \bx\in[0,1]^d,\\
	\displaystyle\phi(\bx,0) = -\frac{1}{2\pi}\cos\left(2\pi\sum_{m=1}^d x_m\right)
	\end{cases}
	\end{equation}
	with periodic boundary conditions.
\end{exa}  

\begin{figure}[h!]
	\centering
%		\subfigure[]{\includegraphics[height=50mm]{cos_p2_n6_t006_surf.pdf}}
	\subfigure[]{\includegraphics[height=50mm]{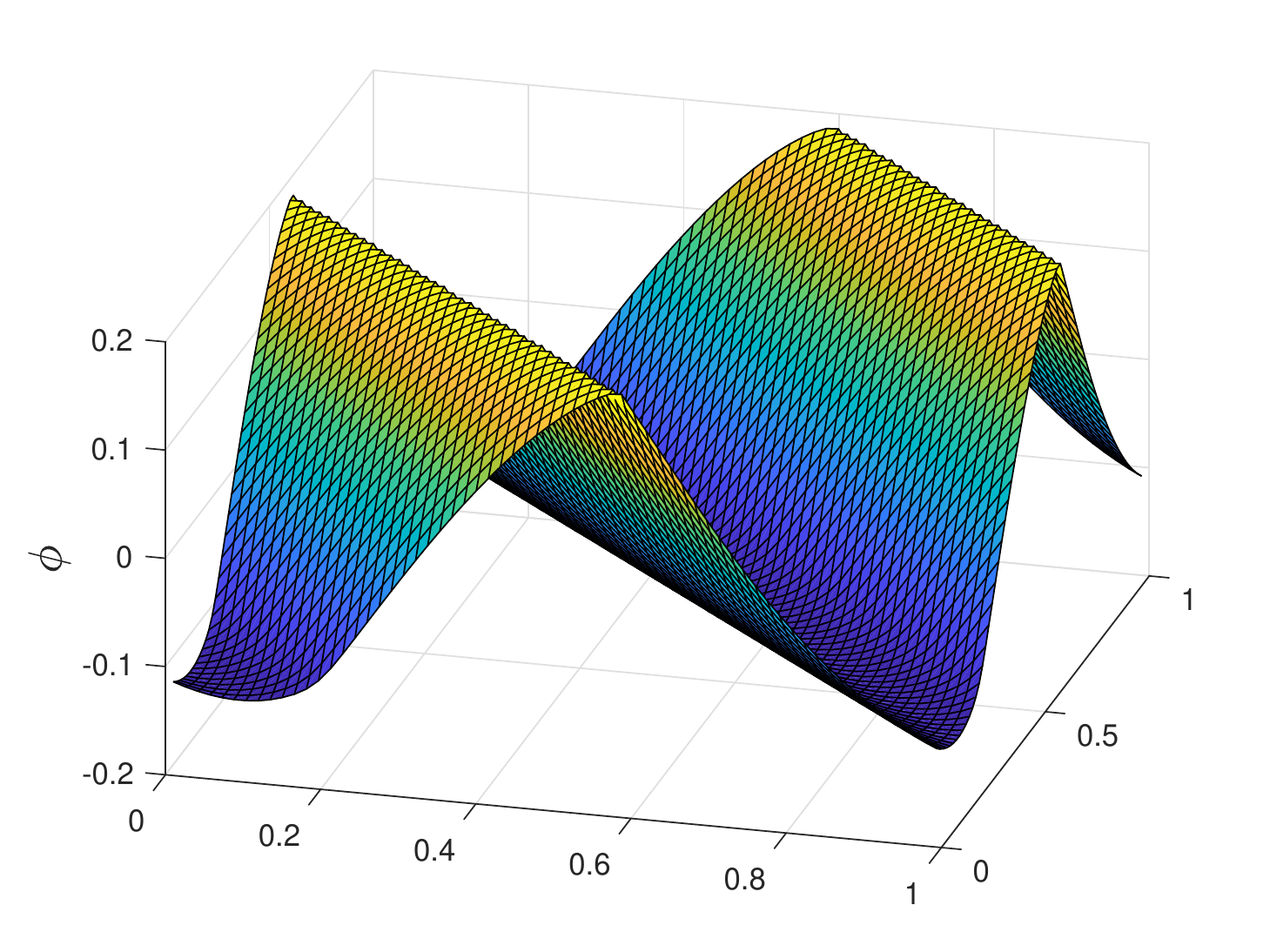}}
	\subfigure[]{\includegraphics[height=50mm]{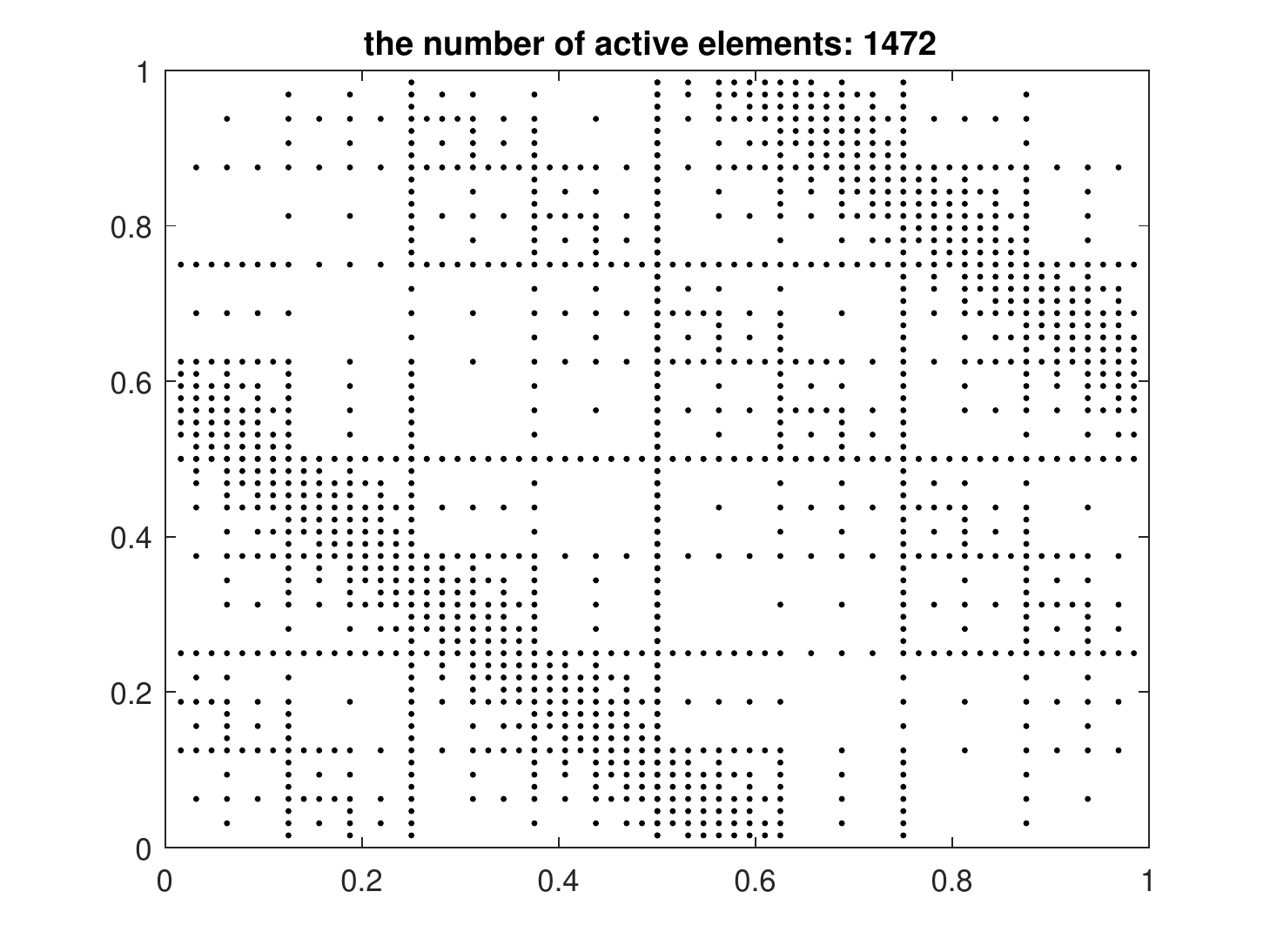}}
	\caption{Example \ref{ex:cos}, $d=2$. $T=0.06$. $k=2$, $M=k$. $N=6$. $\epsilon$=$10^{-5}$.  (a) Numerical solutions on adaptive grids. (b) Active elements.  \label{fig:cos}}
\end{figure}

\begin{figure}[h!]
	\centering
	\subfigure[]{\includegraphics[height=50mm]{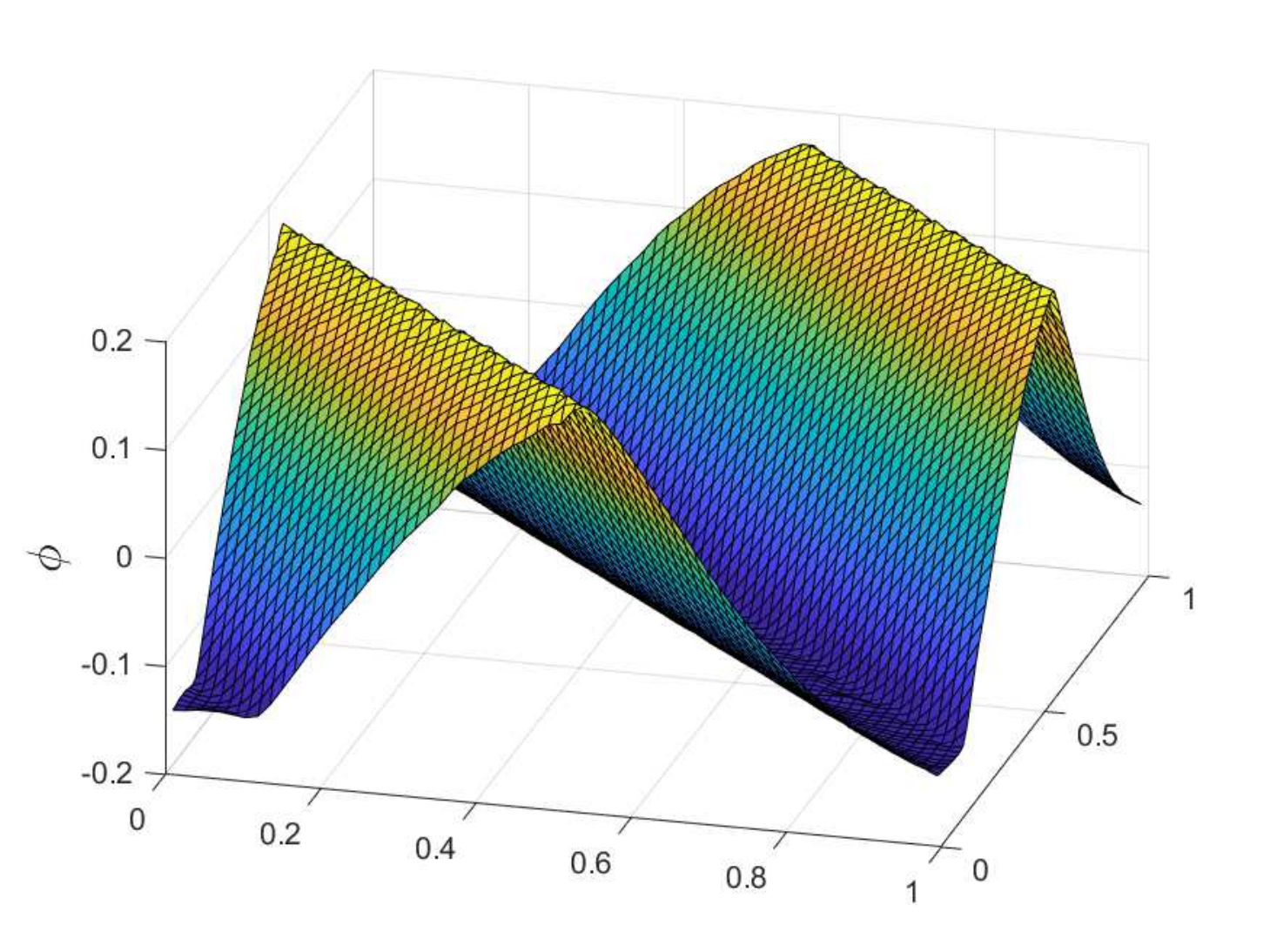}}
	\subfigure[]{\includegraphics[height=50mm]{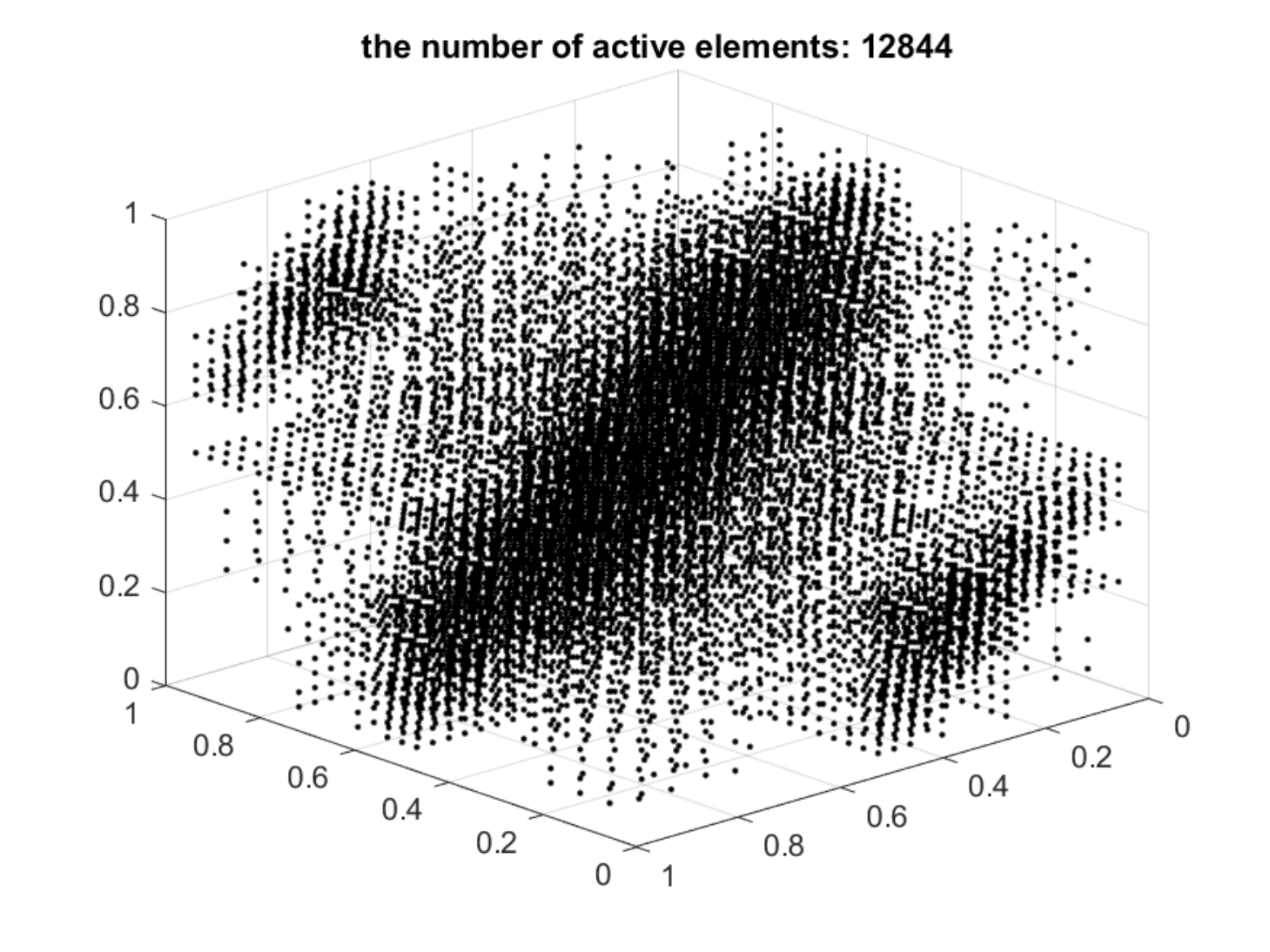}}
	\caption{Example \ref{ex:cos}, $d=3$. $k=2$, $M=3$. $T=0.03$. $N=6$. $\epsilon$=$10^{-6}$. (a) 2D-cuts of the numerical solution at $x_3=0$. (b) Active elements.  \label{fig:cos_3d}}
\end{figure}

In Table \ref{tb:cos_adaptive_d2}-\ref{tb:cos_adaptive_d3}, we report the convergence rates for the adaptive method for $d=2$ and $d=3$ at $T=0.01$ and $T=0.005$, respectively. Similar results are observed to the previous example. In Figure \ref{fig:cos}, we report the solution profile together with the active elements used at $T=0.06$ when the viscosity solution has developed nonsmooth corners. In this simulations, we set $N=6$ and $\epsilon=10^{-5}$. Again, the adaptive method is able to efficiently and correctly capture the sharp corners. In Figure \ref{fig:cos_3d}, we plot the results for $d=3$ at $T=0.03$ with configuration parameters $k=2$, $M=3$, maximum level $N=6$, and $\epsilon=10^{-6}$. High resolution result is observed.

\begin{table}[!hbp]
	\centering
	\caption{Example \ref{ex:cos}, $d=2$. Adaptive sparse grid. $T=0.01$. $M=k$.}
	\label{tb:cos_adaptive_d2}
	\begin{tabular}{c|c|c|c|c|c}
		\hline
		& $\epsilon$ & DoF & L$^2$-error  & $R_{\epsilon}$ & $R_{\textrm{DoF}}$\\
		\hline
		\multirow{5}{3em}{$k=1$}
&	1.00E-03	&	464	&	1.47E-03	&		&		\\
&	1.00E-04	&	1616	&	4.60E-04	&	0.51	&	0.93	\\
&	1.00E-05	&	3840	&	1.66E-04	&	0.44	&	1.18	\\
&	1.00E-06	&	9056	&	2.37E-05	&	0.85	&	2.27	\\
&	1.00E-07	&	17440	&	7.86E-06	&	0.48	&	1.68	\\
		
		\hline
		\hline
		\multirow{5}{3em}{$k=2$}
&	1.00E-03	&	288	&	1.43E-03	&		&		\\
&	1.00E-04	&	720	&	3.20E-04	&	0.65	&	1.64	\\
&	1.00E-05	&	1656	&	9.24E-05	&	0.54	&	1.49	\\
&	1.00E-06	&	3924	&	1.79E-05	&	0.71	&	1.90	\\
&	1.00E-07	&	8406	&	4.00E-06	&	0.65	&	1.97	\\

		\hline
		\hline
		\multirow{5}{3em}{$k=3$}
&	1.00E-03	&	192	&	1.53E-03	&		&		\\
&	1.00E-04	&	512	&	1.85E-04	&	0.92	&	2.15	\\
&	1.00E-05	&	960	&	4.56E-05	&	0.61	&	2.23	\\
&	1.00E-06	&	2048	&	1.50E-05	&	0.48	&	1.47	\\
&	1.00E-07	&	3968	&	1.64E-06	&	0.96	&	3.34	\\
		\hline
	\end{tabular}
\end{table}

\begin{table}[!hbp]
	\centering
	\caption{Example \ref{ex:cos}, $d=3$. Adaptive sparse grid. $T=0.005$. $M=k$.}
	\label{tb:cos_adaptive_d3}
	\begin{tabular}{c|c|c|c|c|c}
		\hline
		& $\epsilon$ & DoF & L$^2$-error  & $R_{\epsilon}$ & $R_{\textrm{DoF}}$\\
		\hline
		\multirow{5}{3em}{$k=1$}
&1.00E-03	&	2432	&	5.22E-03	&	--	&	--	\\
&1.00E-04	&	15680	&	2.70E-03	&	0.29	&	0.35	\\
&1.00E-05	&	46768	&	1.11E-03	&	0.39	&	0.82	\\
&1.00E-06	&	151480	&	4.62E-04	&	0.38	&	0.74	\\
&1.00E-07	&	391008	&	1.66E-04	&	0.45	&	1.08	\\

		\hline
		\hline
		\multirow{5}{3em}{$k=2$}
&1.00E-03	&	2646	&	5.61E-03	&	--	&	--	\\
&1.00E-04	&	8127	&	1.87E-03	&	0.48	&	0.98	\\
&1.00E-05	&	21492	&	6.16E-04	&	0.48	&	1.15	\\
&1.00E-06	&	55296	&	1.24E-04	&	0.70	&	1.70	\\
&1.00E-07	&	154926	&	3.75E-05	&	0.52	&	1.16	\\

		\hline
		\hline
		\multirow{5}{3em}{$k=3$}
&1.00E-03	&	1664	&	2.87E-03	&	--	&	--	\\
&1.00E-04	&	6400	&	1.09E-03	&	0.42	&	0.72	\\
&1.00E-05	&	17920	&	1.60E-04	&	0.83	&	1.87	\\
&1.00E-06	&	44416	&	4.19E-05	&	0.58	&	1.47	\\
&1.00E-07	&	186368	&	5.67E-06	&	0.87	&	1.40	\\
\hline
	\end{tabular}
\end{table}

\begin{figure}[h!]
	\centering
	%\subfigure[]{\includegraphics[height=50mm]{nonlinear_p2_n6_t02_surf.pdf}}
	\subfigure[]{\includegraphics[height=50mm]{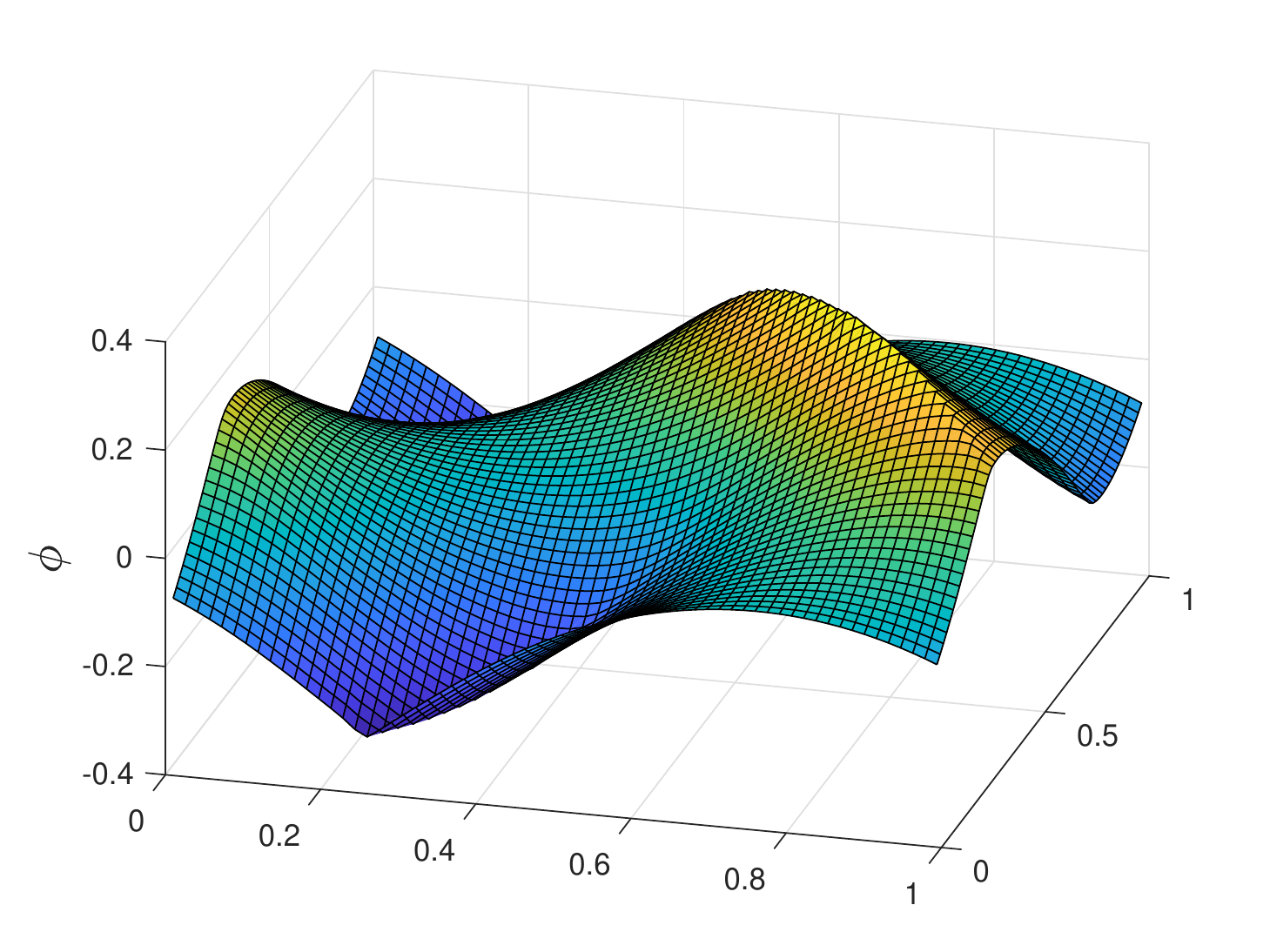}}
	\subfigure[]{\includegraphics[height=50mm]{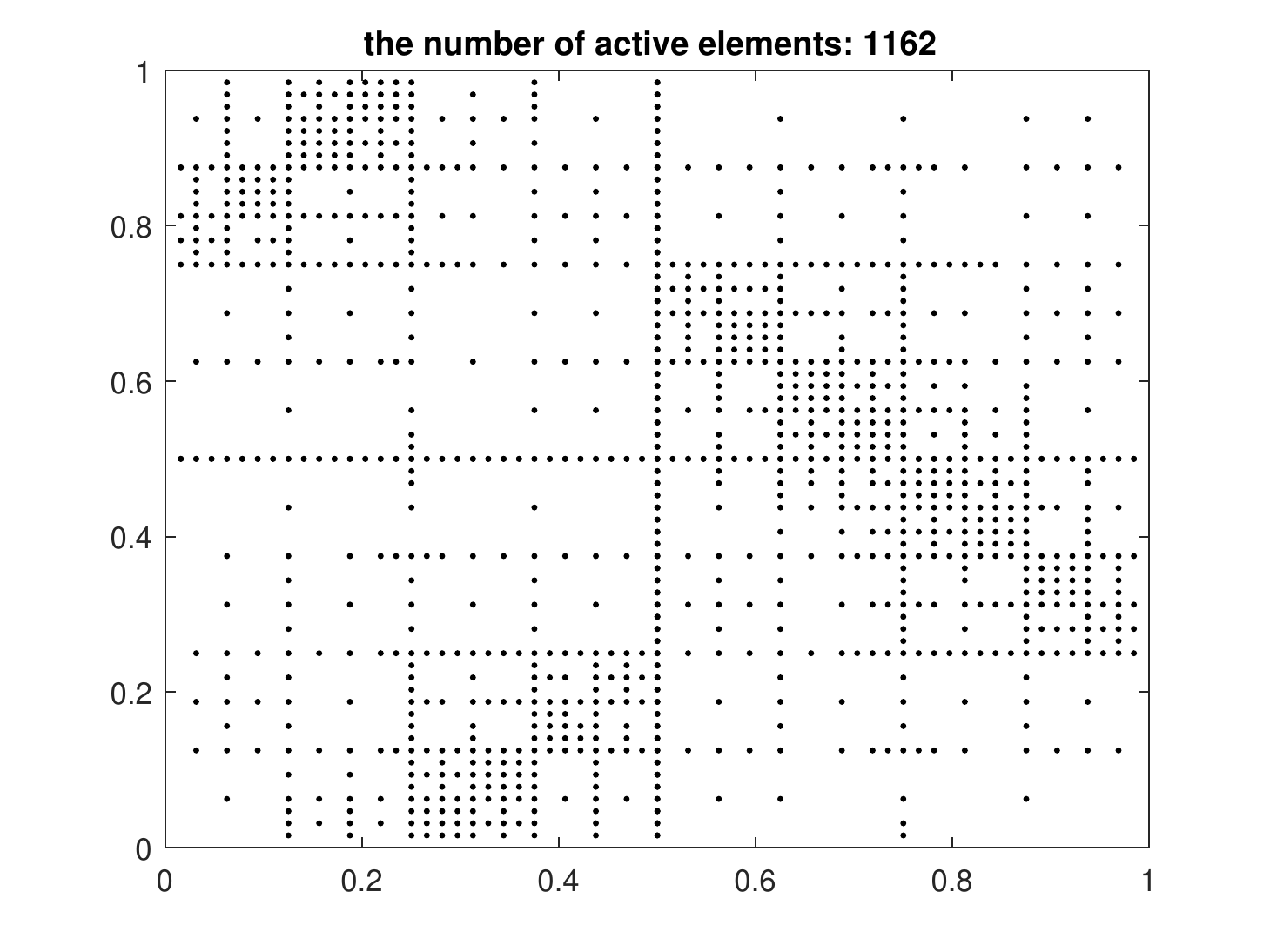}}
	\caption{Example \ref{ex:nonlinear}. $T=0.2$. $N=6$. $k=2$, $M=3$. $\epsilon$=$10^{-5}$. (a) Numerical solution profile. (b) Active elements.  \label{fig:nonlinear}}
\end{figure}

\begin{exa}\label{ex:nonlinear} We consider the following two-dimensional nonlinear problem
	\begin{equation}
	\begin{cases}
	\phi_t + \phi_{x_1}\phi_{x_2} = 0,\quad \bx\in[0,1]^2, \\
	\displaystyle\phi(\bx,0) = -\frac{1}{2\pi}\left(\sin(2\pi x_1)+\cos(2\pi x_2)\right)
	\end{cases}
	\end{equation}
	with periodic boundary conditions.
\end{exa}
Note that unlike the previous two examples, the problem is genuinely nonlinear, and the Hamiltonian is smooth but nonconvex. When $T=0.03$, the solution is still smooth, and we are able to test the convergence for the adaptive method. In the simulation, we set maximum level $N=6$, $k=2$, $M=3$. It is observed in Table \ref{tb:nonlinear_adap} that the method is able to achieve very accurate results by using a few DoFs. The convergence performance is similar to the previous examples. In Figure \ref{fig:nonlinear}, we plot the solution at $T=0.2$, when the viscosity solution becomes nonsmooth. It is observed that the adaptive method captures the corners correctly and efficiently, as compared with the results by other popular methods, see e.g. \cite{li2010central,cheng2014hj}.

\begin{table}[!hbp]
	\centering
	\caption{Example \ref{ex:nonlinear}.   Adaptive sparse grid. $T=0.03$. $M=k$.}
	\label{tb:nonlinear_adap}
	\begin{tabular}{c|c|c|c|c|c}
		\hline
		& $\epsilon$ & DoF & L$^2$-error  & $R_{\epsilon}$ & $R_{\textrm{DoF}}$\\
		\hline
		\multirow{5}{3em}{$k=1$}
		&	1.00E-03	&	180	&	9.78E-03	&	--	&	--	\\
		&	1.00E-04	&	448	&	2.05E-03	&	0.68	&	1.71	\\
		&	1.00E-05	&	952	&	1.26E-03	&	0.21	&	0.64	\\
		&	1.00E-06	&	1296	&	2.24E-04	&	0.75	&	5.60	\\
		&	1.00E-07	&	2952	&	2.74E-05	&	0.91	&	2.56	\\

		\hline
		\hline
		\multirow{5}{3em}{$k=2$}
		&	1.00E-03	&	135	&	2.73E-03	&	--	&	--	\\
		&	1.00E-04	&	306	&	3.95E-04	&	0.84	&	2.36	\\
		&	1.00E-05	&	594	&	1.97E-04	&	0.30	&	1.05	\\
		&	1.00E-06	&	1224	&	4.41E-05	&	0.65	&	2.07	\\
		&	1.00E-07	&	2565	&	1.31E-05	&	0.53	&	1.64	\\

		\hline
		\hline
		\multirow{5}{3em}{$k=3$}
		&	1.00E-03	&	112	&	5.38E-04	&	--	&	--	\\
		&	1.00E-04	&	256	&	1.67E-04	&	0.51	&	1.42	\\
		&	1.00E-05	&	560	&	4.29E-05	&	0.59	&	1.73	\\
		&	1.00E-06	&	832	&	1.21E-05	&	0.55	&	3.20	\\
		&	1.00E-07	&	1280	&	1.25E-06	&	0.99	&	5.27	\\
		
		\hline
	\end{tabular}
\end{table}

\begin{figure}[h!]
	\centering
	\subfigure[]{\includegraphics[height=50mm]{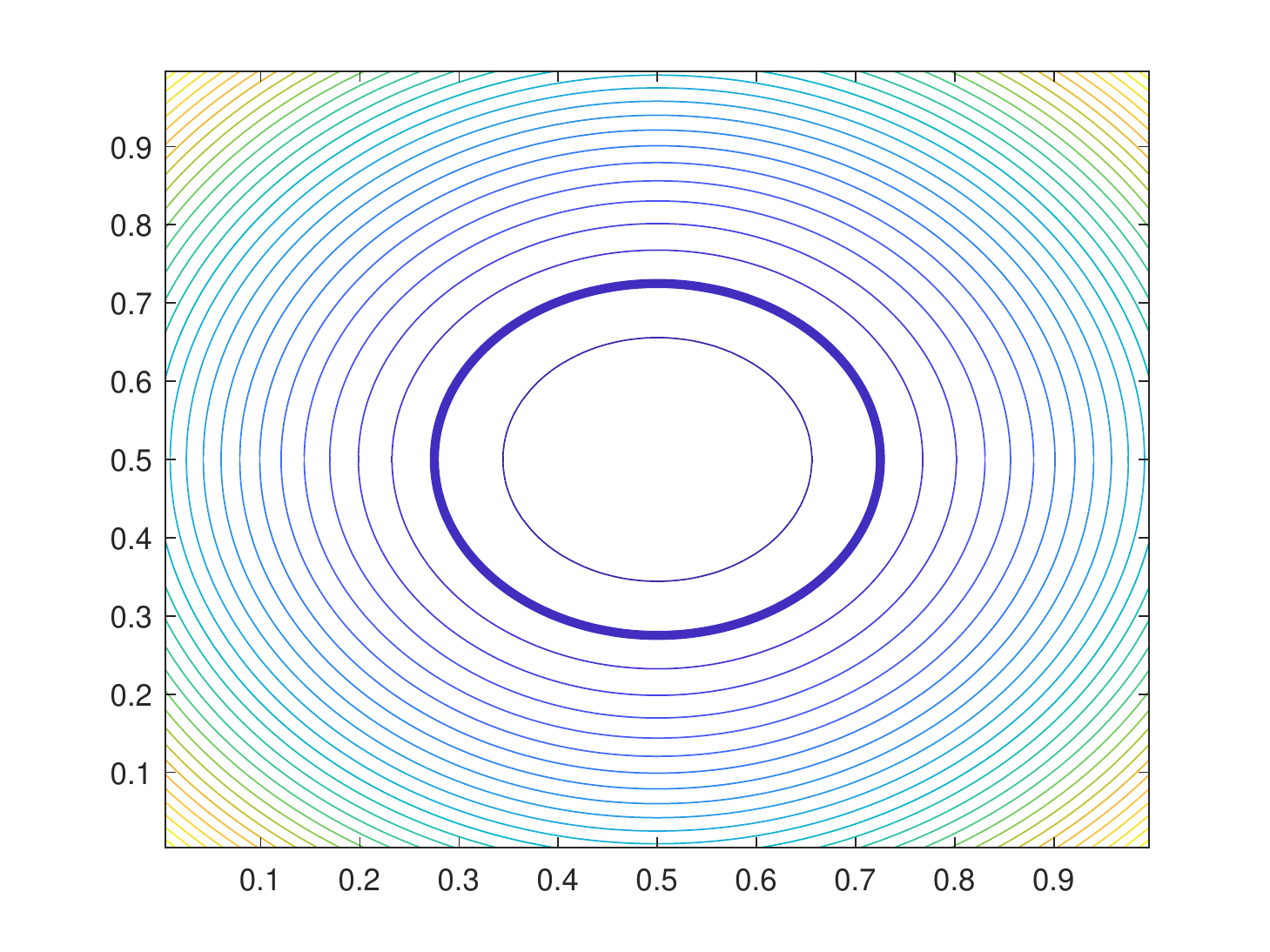}}
	\subfigure[]{\includegraphics[height=50mm]{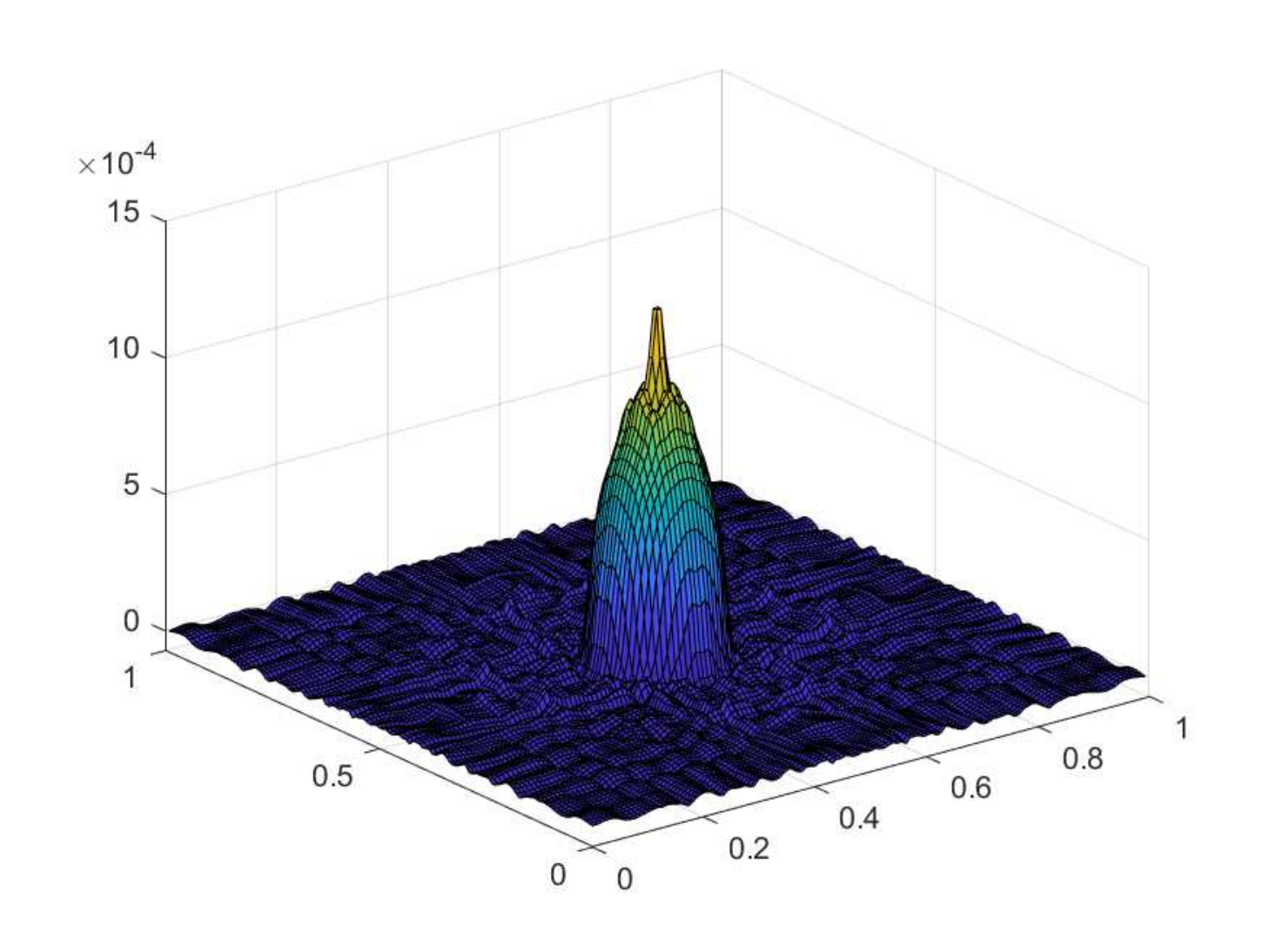}}
	\subfigure[]{\includegraphics[height=50mm]{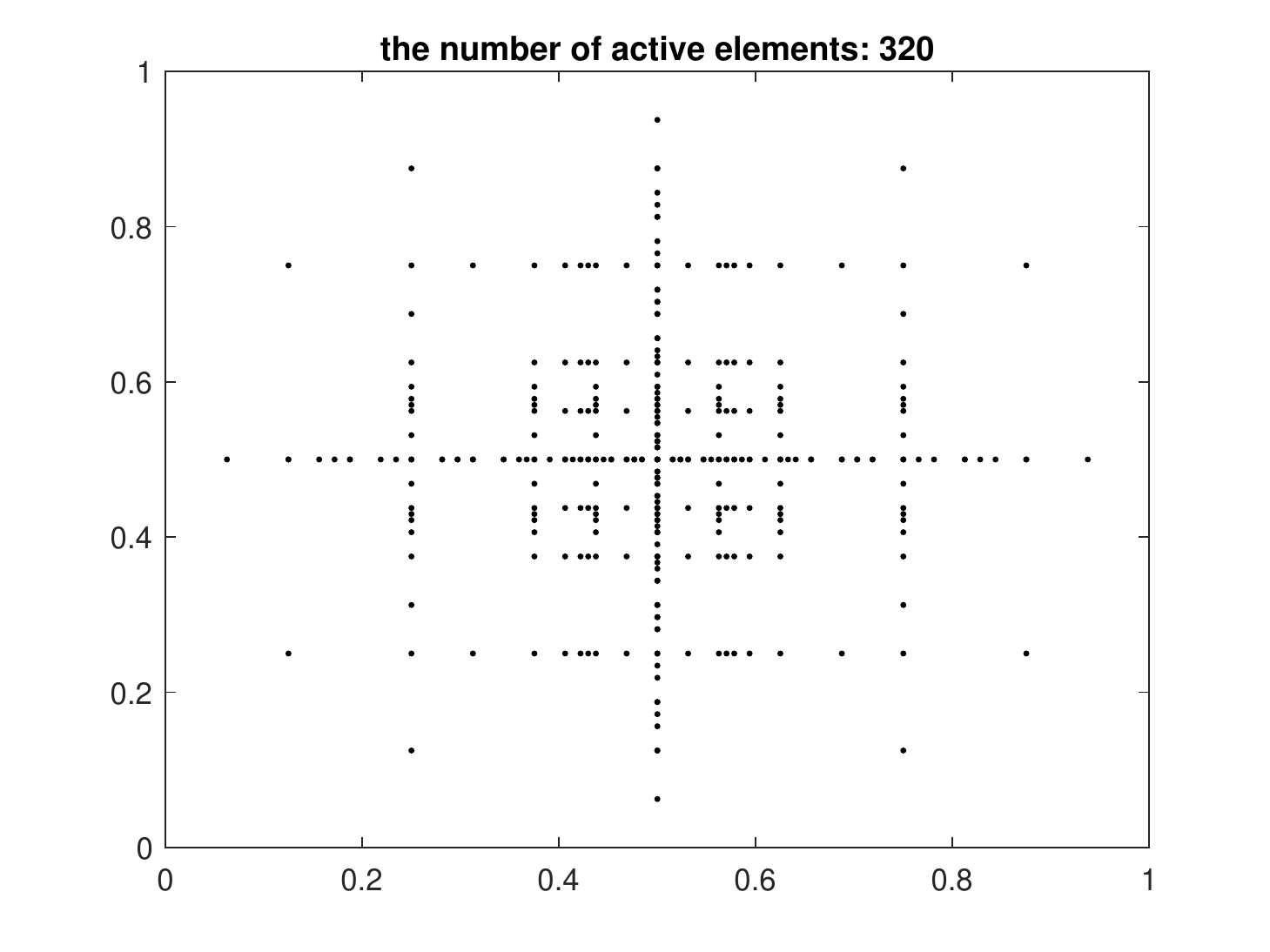}}
	\caption{Example \ref{ex:eikonal}. $T=0.1$. $k=2$, $M=3$. $N=7$. $\epsilon$=$10^{-7}$. (a) Contour plot of the numerical solution. (b) Numerical error distribution. (c) Active elements.  \label{fig:eikonal}}
\end{figure}

\begin{exa}\label{ex:eikonal}We consider the classic nonlinear Eikonal equation
	\begin{equation}
	\left\{\begin{array}{l}\displaystyle
	\phi_{t} + \|\nabla\phi\|=0 ,\quad \bx\in[0,1]^d\\
	\phi(\bx, 0)=g(\|\bx-\ba\|)
	\end{array}\right.
	\end{equation}
	where $\ba=(0.5,0.5,\ldots,0.5)$ and 
	$$g(z) = \frac{1}{2r_0}(z^2-r_0^2),\quad r_0 = \frac18.$$
	An outflow boundary condition is imposed. 
	The viscosity solution is 
	$$
	\phi(\bx,t) = g\left(\max\left(\|\bx\| -t,0\right)\right),
	$$
	which is clearly $C^1$ smooth. 
\end{exa}
 One additional challenge of this problem is that the Hamiltonian is not smooth, making the DG formulation unstable if the numerical quadrature is not sufficiently accurate, as mentioned in previous section. To circumvent the difficulty, we propose to employ a regularized Hamiltonian as follows.
\begin{equation}
\label{eikonal_regu}
\tilde H(\nabla\phi) = \begin{cases}
\|\nabla\phi\|, &\text{if } \|\nabla\phi\| \ge \delta  \\
\frac{1}{2\delta} \|\nabla\phi\|^2 + \frac12\delta, &\text{otherwise}.
\end{cases}
\end{equation}
It can be easily verified that $\tilde H$ is $C^1$. In the simulation, we choose $\delta = 2h$, where $h$ is the mesh size, hence the regularization will not affect the accuracy of the original method. We employ the regularized Hamiltonian for all tests, while we notice that it is only required for $k>1$.
In Table \ref{tb:eikonal_sparse_d2}, we summarize the convergence study for the adaptive method for $d=2,\,3,\,4$ and $k=1,\,2$.   It is observed that the convergence rates $R_{\textrm{DoF}}$ and $R_{\epsilon}$ are similar for $k=1$ and $k=2$, which is unsurprising, since the viscosity solution is only $C^1$. Meanwhile, the error magnitude by $k=2$ is still much smaller than that by $k=1$ with the same number of DoF, demonstrating the efficiency of method with high  order accuracy. In Figure \ref{fig:eikonal}, we report the contour plot of the numerical solution with $N=7$, $k=2$, $M=3$, $\epsilon = 10^{-7}$, $d=2$. We observe that the rarefaction wave developed at the center of domain is correctly captured by the adaptive method.
 We also highlight the level set of $\phi=0$.

 \begin{table}[!hbp]
 	\centering
 	\caption{Example \ref{ex:eikonal}, $d=2,\,3,\,4$. Adaptive sparse grid. $T=0.1$. $M=k+1$.}
 	\label{tb:eikonal_sparse_d2}
 	\begin{tabular}{c|c|c|c|c|c|c|c|c|c}
 		\hline
  \multicolumn{2}{c|}{} &  \multicolumn{4}{c|}{$k=1$} &\multicolumn{4}{c}{$k=2$}\\
 			\hline
 		& $\epsilon$ & DoF & L$^2$-error  & $R_{\epsilon}$ & $R_{\textrm{DoF}}$& DoF & L$^2$-error  & $R_{\epsilon}$ & $R_{\textrm{DoF}}$\\
 		\hline
 		\multirow{5}{3em}{$d=2$}
&1.00E-03	&	236	&	2.25E-02	&		&		&	72	&	5.42E-03	&		&	\\
&1.00E-04	&	496	&	5.39E-03	&	0.62	&	1.92	&	108	&	3.59E-03	&	0.18	&	1.01		\\
&1.00E-05	&	1056	&	2.93E-03	&	0.26	&	0.81	&	324	&	1.09E-03	&	0.52	&	1.09		\\
&1.00E-06	&	1904	&	9.27E-04	&	0.50	&	1.95	&	900	&	4.41E-04	&	0.39	&	0.89		\\
&1.00E-07	&	5496	&	2.43E-04	&	0.58	&	1.26	&	2880	&	1.23E-04	&	0.56	&	1.10		\\
%&1.00E-08	&	12212	&	2.57E-04	&	-0.02	&	-0.07	&	8928	&	8.99E-05	&	0.13	&	0.27		\\

 		\hline
 		\hline
 		\multirow{5}{3em}{$d=3$}
&1.00E-03	&	680	&	2.31E-02	&		&		&	108	&	6.46E-03	&		&			\\
&1.00E-04	&	1472	&	7.64E-03	&	0.48	&	1.43	&	351	&	3.18E-03	&	0.31	&	0.60		\\
&1.00E-05	&	2968	&	4.12E-03	&	0.27	&	0.88	&	1026	&	1.53E-03	&	0.32	&	0.68		\\
&1.00E-06	&	5080	&	1.72E-03	&	0.38	&	1.63	&	2970	&	5.76E-04	&	0.43	&	0.92		\\
&1.00E-07	&	23272	&	5.00E-04	&	0.54	&	0.81	&	11610	&	2.26E-04	&	0.41	&	0.69		\\

 		\hline
 		\hline
 		\multirow{5}{3em}{$d=4$} 
 		&	1.00E-03	&	1872	&	2.30E-02	&		&		&	405	&	4.61E-03	&		&		\\
 		&	1.00E-04	&	3792	&	2.52E-02	&	-0.04	&	-0.13	&
 			1053	&	2.90E-03	&	0.20	&	0.48	\\
 		&	1.00E-05	&	8944	&	1.25E-02	&	0.30	&	0.81	&
 			3159	&	1.16E-03	&	0.40	&	0.84	\\
 		&	1.00E-06	&	10624	&	2.59E-03	&	0.68	&	9.15	&
 			12312	&	5.71E-04	&	0.31	&	0.52	\\
 		&	1.00E-07	&	-	&	-	&	-	& -	&	55080	&	2.03E-04	&	0.45	&	0.69
 		\\

 				\hline
 	\end{tabular}
 \end{table}

\begin{figure}[h!]
	\centering
	\subfigure[]{\includegraphics[height=50mm]{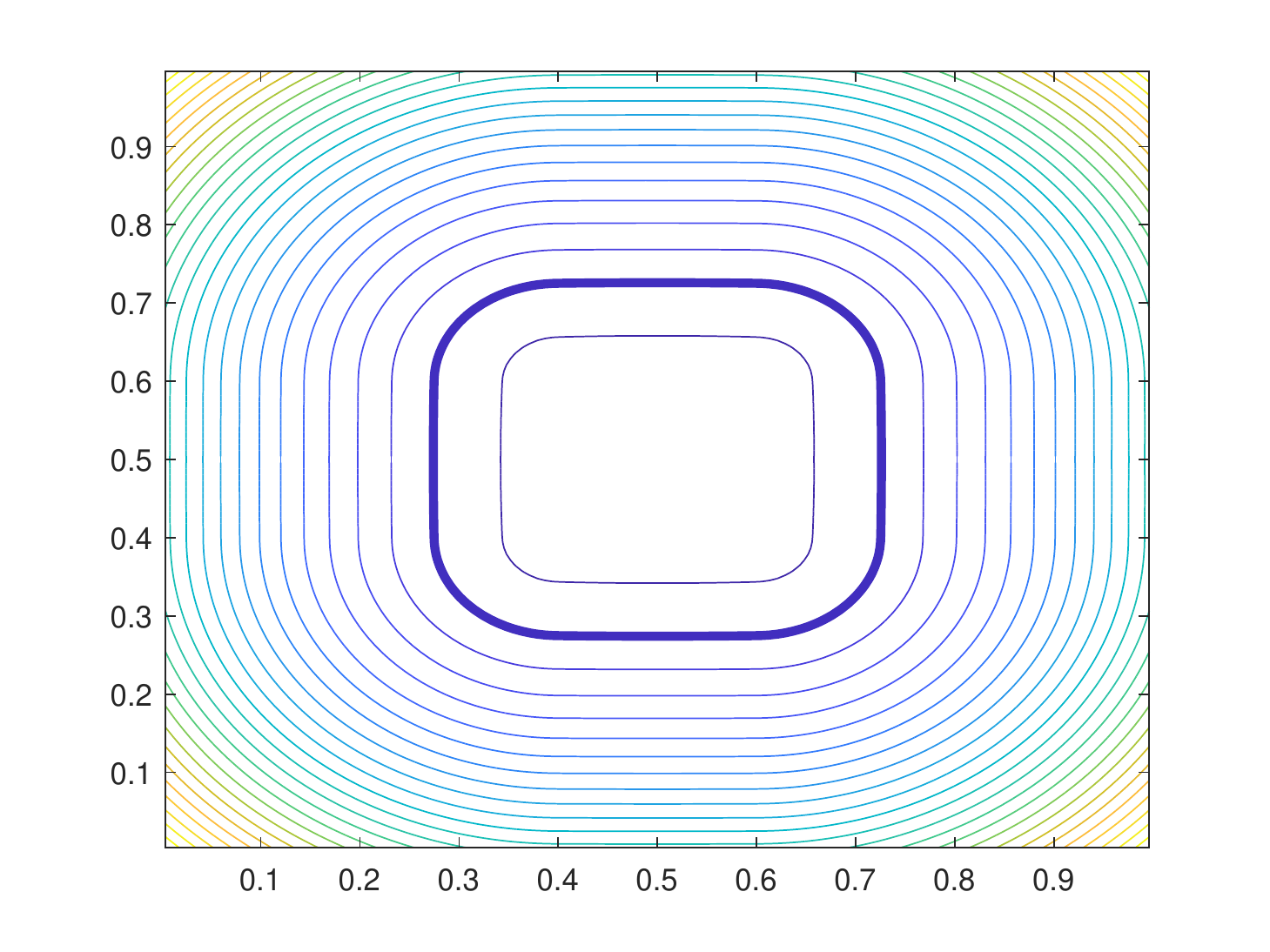}}
	\subfigure[]{\includegraphics[height=50mm]{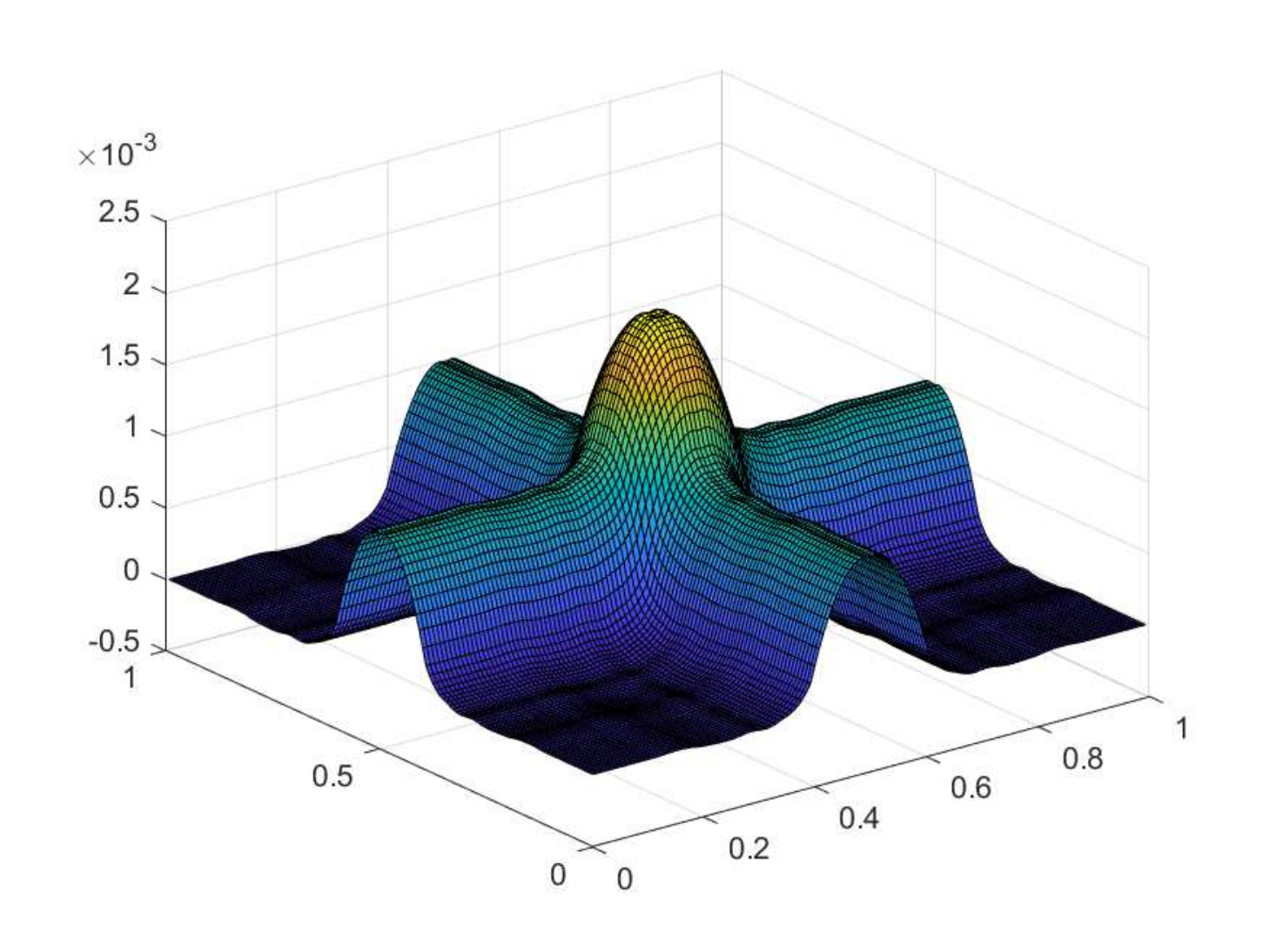}}
	\subfigure[]{\includegraphics[height=50mm]{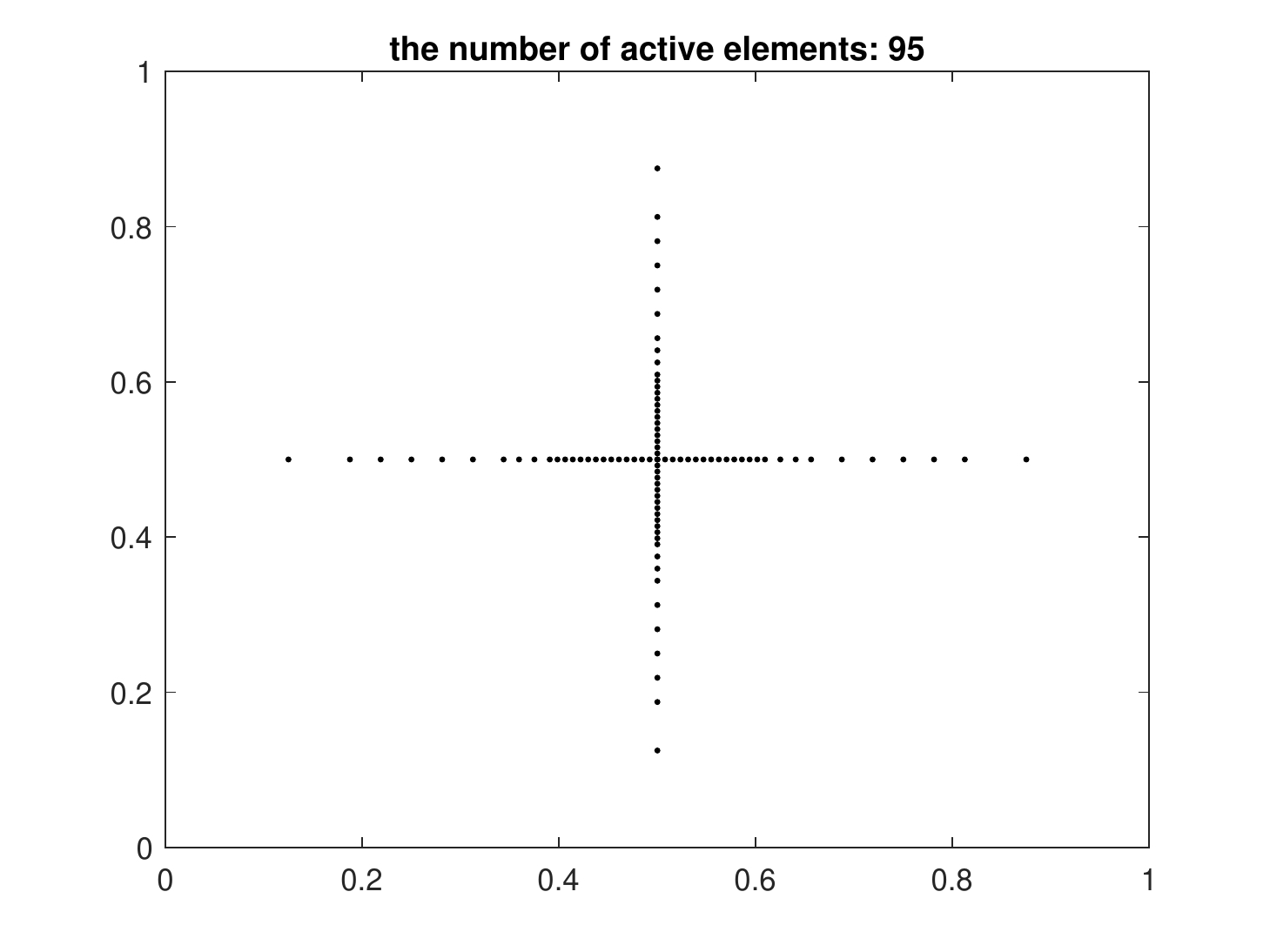}}
	\caption{Example \ref{ex:HJB}. $T=0.1$. $k=2$, $M=4$. $N=7$. $\epsilon$=$10^{-7}$. (a) Contour plot of the numerical solution. (b) Numerical error distribution. (c) Active elements.  \label{fig:hjb}}
\end{figure}

\begin{exa}\label{ex:HJB}
	In this example, we consider the following HJB equation \cite{bokanowski2013adaptive}
\begin{equation}\label{ex:HJB_ex}
\left\{\begin{array}{l}\displaystyle
\phi_t +\max_{\bb\in\mathcal{B}}\left(\sum_{m=1}^d b_m\cdot\nabla\phi\right) = 0,\quad \bx\in[0,1]^d,\\
\phi(\bx, 0)=g(\|\bx-\ba\|),
\end{array}\right.
\end{equation}	
	where $\ba=(0.5,0.5,\ldots,0.5)$ and $\mathcal{B} = \{\bb=(b_1,b_2,\ldots,b_d),\,b_m=\pm1\} $ is a set of $2^d$ vectors corresponding to $2^d$ possible controls. 
%	$P$ denotes the $d\times d$ matrix as follows
%	$$P:=\left[\begin{array}{cccccc}
%1 & 0 & 0 & \ldots & \ldots & 0 \\
%-1 & -1 & 0 & \ldots & \ldots & 0 \\
%0 & 0 & 1 & & & \\
%\vdots & & 0 & \ddots & & \\
%\vdots & & \vdots & & \ddots & \\
%0 & \ldots & 0 & & & 1
%\end{array}\right],$$
%and $P_m$ denotes the $m$-th column of $P$. 
The function $g(z)$ is the same as in the example \ref{ex:eikonal}. Note that this HJB equation is equivalent to the following HJ equation
\begin{equation}
\label{ex:HJ_ex}
\left\{\begin{array}{l}\displaystyle
\phi_t + \sum_{m=1}^d|\phi_{x_m}| = 0,\quad \bx\in[0,1]^d\\
\phi(\bx, 0)=g(\|\bx-\ba\|),
\end{array}\right.
\end{equation}	
The exact solution can be hence derived from \eqref{ex:HJ_ex}:
$$
\phi(\bx, t)=g(\|(\bx-\ba))_t^\star\|).
$$ 
Here, for a vector $\mathbf{c}$, $\mathbf{c}_t^\star: = \min (\max (0, \mathbf{c}-t), \mathbf{c}+t)$ in the component-wise sense.
We apply the adaptive algorithm to simulate \eqref{ex:HJ_ex}. The outflow boundary conditions are imposed. Note that the Hamiltonian is nonsmooth as with the Eikonal equation, and hence we regularize the absolute function using the technique \eqref{eikonal_regu} to ensure stability.
In Figure \ref{fig:hjb}, we plot the solution with configuration $k=2$, $M=4$, $N=7$, $\epsilon=10^{-7}$. Note that the viscosity solution is $C^1$, a rarefaction wave opens up at the center of the domain, which is well captured by the method. In Table \ref{tb:hjb_sparse_d2}, we summarize the convergence study for $d=2,\,3,\,4$ and $k=1,\,2$. Note that when $\epsilon = 10^{-7}$, the error does not decay anymore, since it has  saturated already with the maximum level $N=7$.

\begin{table}[!hbp]
	\centering
	\caption{Example \ref{ex:HJB}, $d=2,\,3,\,4$. Adaptive sparse grid. $T=0.1$. $M=k+2$.}
	\label{tb:hjb_sparse_d2}
	\begin{tabular}{c|c|c|c|c|c|c|c|c|c}
		\hline
		\multicolumn{2}{c|}{} &  \multicolumn{4}{c|}{$k=1$} &\multicolumn{4}{c}{$k=2$}\\
		\hline
		& $\epsilon$ & DoF & L$^2$-error  & $R_{\epsilon}$ & $R_{\textrm{DoF}}$& DoF & L$^2$-error  & $R_{\epsilon}$ & $R_{\textrm{DoF}}$\\
		\hline
		\multirow{6}{3em}{$d=2$}
&	1.00E-03	&	204	&	4.17E-03	&		&		&	63	&	8.62E-03	&		&		\\
&	1.00E-04	&	444	&	1.62E-03	&	0.41	&	1.21	&	135	&	1.89E-03	&	0.66	&	1.99	\\
&	1.00E-05	&	860	&	7.01E-04	&	0.36	&	1.27	&	207	&	6.26E-04	&	0.48	&	2.59	\\
&	1.00E-06	&	876	&	6.43E-04	&	0.04	&	4.69	&	459	&	4.24E-04	&	0.17	&	0.49	\\
&	1.00E-07	&	924	&	6.58E-04	&	-0.01	&	-0.43	&	855	&	3.97E-04	&	0.03	&	0.11	\\

		\hline
		\hline
		\multirow{5}{3em}{$d=3$}
&	1.00E-03	&	608	&	5.44E-03	&		&		&	270	&	1.12E-02	&		&		\\
&	1.00E-04	&	1328	&	2.12E-03	&	0.41	&	1.20	&	594	&	2.98E-03	&	0.58	&	1.68	\\
&	1.00E-05	&	2576	&	9.15E-04	&	0.37	&	1.27	&	918	&	7.89E-04	&	0.58	&	3.05	\\
&	1.00E-06	&	2624	&	8.39E-04	&	0.04	&	4.74	&	2052	&	4.36E-04	&	0.26	&	0.74	\\
&	1.00E-07	&	2768	&	8.56E-04	&	-0.01	&	-0.38	&	3510	&	3.93E-04	&	0.05	&	0.19	\\

		\hline
		\hline
		\multirow{5}{3em}{$d=4$} 
&	1.00E-03	&	1616	&	6.65E-03	&		&		&	1053	&	1.35E-02	&		&		\\
&	1.00E-04	&	3536	&	2.60E-03	&	0.41	&	1.20	&	1701	&	6.17E-03	&	0.34	&	1.63	\\
&	1.00E-05	&	6864	&	1.12E-03	&	0.37	&	1.27	&	2997	&	1.24E-03	&	0.70	&	2.84	\\
&	1.00E-06	&	6992	&	1.02E-03	&	0.04	&	4.69	&	6237	&	4.90E-04	&	0.40	&	1.26	\\
&	1.00E-07	&	7376	&	1.05E-03	&	-0.01	&	-0.38	&	12069	&	4.20E-04	&	0.07	&	0.23	\\

		\hline
	\end{tabular}
\end{table}

\end{exa}

%%%%%%%%%%%%%%%%%%%%%

\begin{figure}[h!]
	\centering
	\subfigure[]{\includegraphics[height=50mm]{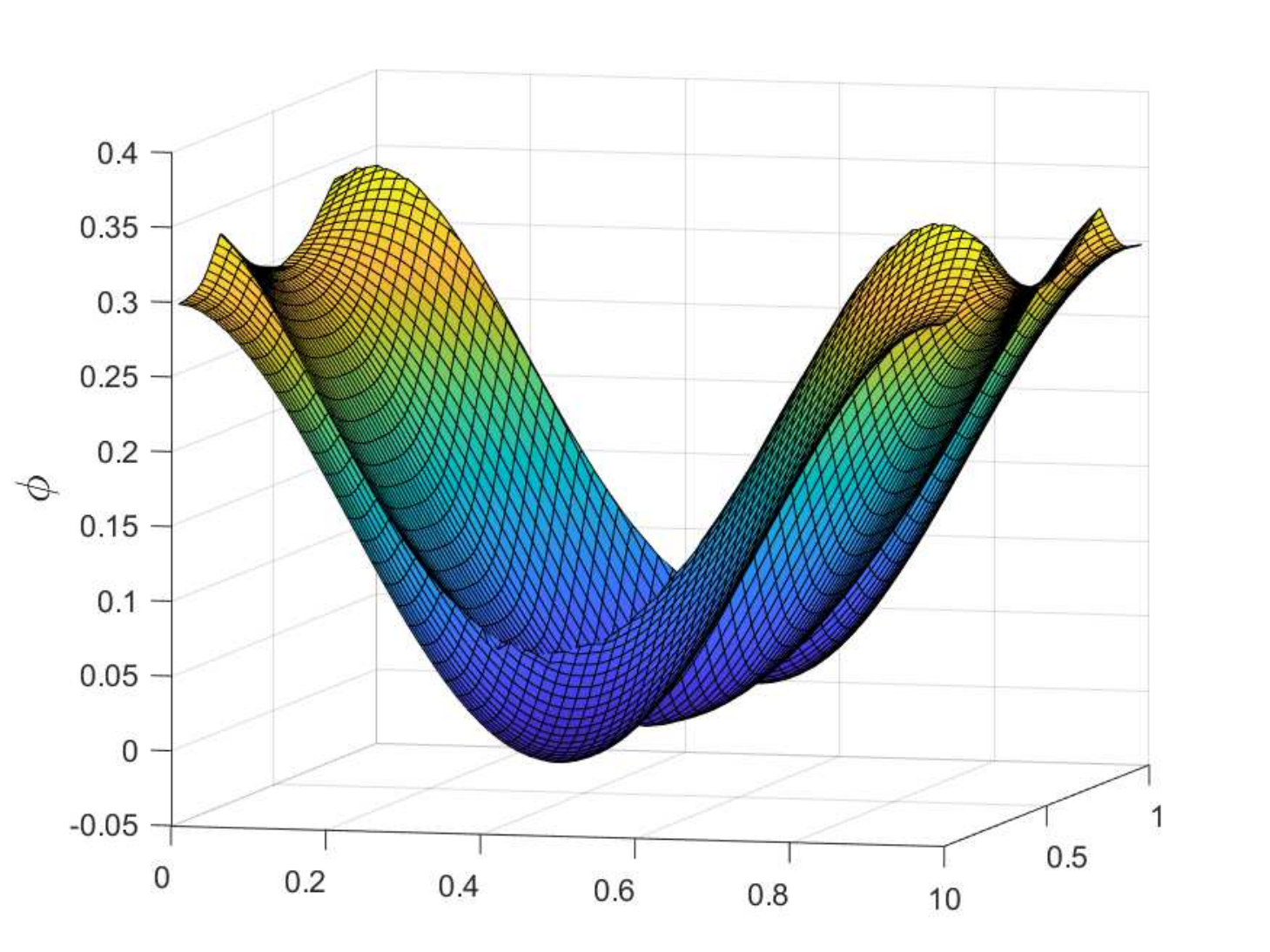}}
	\subfigure[]{\includegraphics[height=50mm]{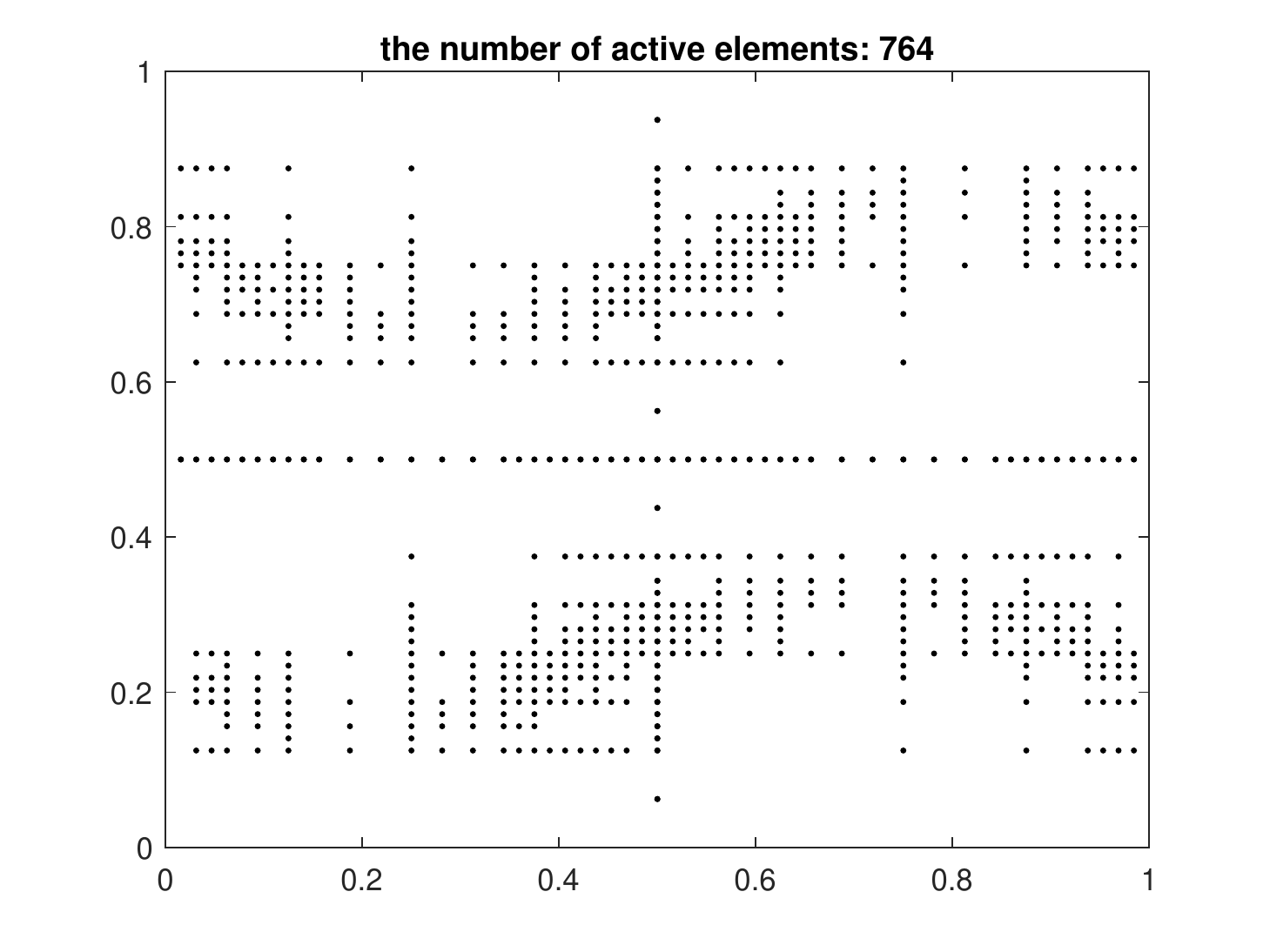}}
     \subfigure[]{\includegraphics[height=50mm]{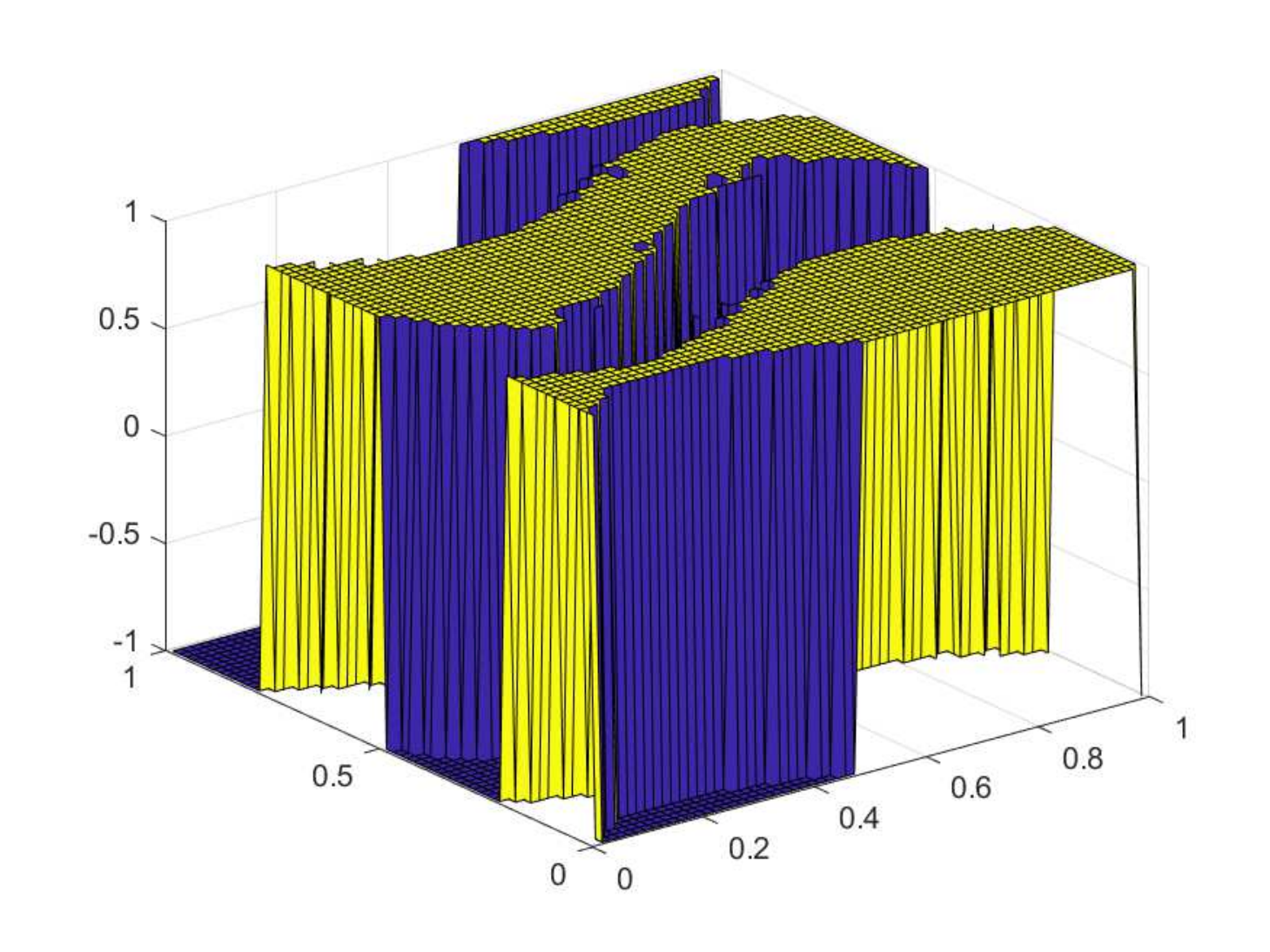}}
	\caption{Example \ref{ex:optimal}, $d=2$. $T=0.15$. $k=2$, $M=4$. $N=6$. $\epsilon$=$10^{-5}$. (a) Numerical solution profile. (b) Active elements. (c) Controls $\text{sign}(\phi_y)$  \label{fig:controls}}
\end{figure}

\begin{exa}\label{ex:optimal}In the last example, we consider the 2D problem related to controlling optimal cost determination \cite{osher1991high}
\begin{equation}
\left\{\begin{array}{l}\displaystyle
\phi_{t}-\sin (2\pi x_2) \phi_{x_1}-\left(\sin (2\pi x_1)+\operatorname{sign}\left(\phi_{x_2}\right)\right) \phi_{x_2}-\frac{1}{2} \sin ^{2}(2\pi x_2) - \cos (2\pi x_1)-1=0 ,\quad \bx\in[0,1]^2\\
\phi(\bx, 0)=0
\end{array}\right.
\end{equation}
\end{exa}
Note that the Hamiltonian is not smooth. In Figure \ref{fig:controls}, we plot the solution profile, the optimal $\text{sign}(\phi_{x_2})$ together with the active elements at final time $T=0.15$. Again, we regularize the Hamiltonian as with previous examples. The adaptive method is able to capture the viscosity solution efficiently, and the numerical results agree with other methods in the literature, e.g. \cite{hu1999discontinuous,li2010central,cheng2014hj,guo2011local,ke2019alternative}.

% section 6
\section{Conclusion}\label{sec:conclusion}

In this work, we proposed an adaptive sparse grid LDG method for solving HJ equations in high dimensions. By incorporating the orthonormal Alpert's multiwavelets as the DG finite element bases, and the interpolatory multiwavelets as efficient multiresolution numerical quadratures, we achieve efficient multiresolution schemes which is suitable for high dimensions. Benchmark numerical tests up to 4D are provided to validate the performance of the method. %Only very few DoFs are needed for high dimension problems. %Future work includes extensions to other kinetic equations in high dimensions. 
%In an effort for promoting reproducible research, 
The code generating the results in this paper can be found at the GitHub link: \url{https://github.com/JuntaoHuang/adaptive-multiresolution-DG}, and it has the capability of computing higher dimensional problems.

% section 7
\section*{Acknowledgment}

We would like to thank Qi Tang and Kai Huang for the assistance and discussion in code implementation. Yingda Cheng would like to thank the support from IPAM to attend the workshop on ``High-dimensional Hamilton-Jacobi PDEs".

\bibliographystyle{abbrv}
%\bibliographystyle{plain}
% \bibliography{bibliography/BigBib,bibliography/adapt_cl,bibliography/ref_cheng,bibliography/ref_cheng_2,bibliography/adaptive,bibliography/some_wave}
% \bibliographystyle{abbrv}
% \bibliography{some_wave,ref_hj,BigBib,adapt_cl,ref_cheng,ref_cheng_2,adaptive}

\begin{thebibliography}{10}

\bibitem{abgrall1996numerical}
R.~Abgrall.
\newblock {Numerical discretization of the first-order Hamilton-Jacobi equation
  on triangular meshes}.
\newblock {\em Communications on Pure and Applied Mathematics},
  49(12):1339--1373, 1996.

\bibitem{alpert1993class}
B.~Alpert.
\newblock A class of bases in {$L$}$^2$ for the sparse representation of
  integral operators.
\newblock {\em SIAM Journal on Mathematical Analysis}, 24(1):246--262, 1993.

\bibitem{bokanowski2013adaptive}
O.~Bokanowski, J.~Garcke, M.~Griebel, and I.~Klompmaker.
\newblock An adaptive sparse grid semi-{L}agrangian scheme for first order
  {H}amilton-{J}acobi {B}ellman equations.
\newblock {\em Journal of Scientific Computing}, 55(3):575--605, 2013.

\bibitem{bungartz2004sparse}
H.-J. Bungartz and M.~Griebel.
\newblock Sparse {G}rids.
\newblock {\em Acta Numerica}, 13:147--269, 2004.

\bibitem{cheng2007discontinuous}
Y.~Cheng and C.-W. Shu.
\newblock {A discontinuous Galerkin finite element method for directly solving
  the Hamilton--Jacobi equations}.
\newblock {\em Journal of Computational Physics}, 223(1):398--415, 2007.

\bibitem{cheng2014hj}
Y.~Cheng and Z.~Wang.
\newblock A new discontinuous {Galerkin} finite element method for directly
  solving the {Hamilton-Jacobi} equations.
\newblock {\em Journal of Computational Physics}, 268:134--153, 2014.

\bibitem{chow2017algorithm}
Y.~Chow, J.~Darbon, S.~Osher, and W.~Yin.
\newblock {Algorithm for overcoming the curse of dimensionality for
  time-dependent non-convex Hamilton--Jacobi equations arising from optimal
  control and differential games problems}.
\newblock {\em Journal of Scientific Computing}, 73(2-3):617--643, 2017.

\bibitem{chow2019algorithm}
Y.~T. Chow, J.~Darbon, S.~Osher, and W.~Yin.
\newblock Algorithm for overcoming the curse of dimensionality for
  state-dependent Hamilton-Jacobi equations.
\newblock {\em Journal of Computational Physics}, 387:376--409, 2019.

\bibitem{DG4}
B.~Cockburn, S.~Hou, and C.-W. Shu.
\newblock The {R}unge-{K}utta local projection discontinuous {G}alerkin finite
  element method for conservation laws. {IV}. {T}he multidimensional case.
\newblock {\em Mathematics of Computation}, 54(190):545--581, 1990.

\bibitem{cockburn1998localDG}
B.~Cockburn and C.-W. Shu.
\newblock The local discontinuous {G}alerkin method for time-dependent
  convection-diffusion systems.
\newblock {\em SIAM Journal on Numerical Analysis}, 35(6):2440--2463, 1998.

\bibitem{crandall1984some}
M.~G. Crandall, L.~C. Evans, and P.-L. Lions.
\newblock {Some properties of viscosity solutions of Hamilton-Jacobi
  equations}.
\newblock {\em Transactions of the American Mathematical Society},
  282(2):487--502, 1984.

\bibitem{crandall1983viscosity}
M.~G. Crandall and P.-L. Lions.
\newblock {Viscosity solutions of Hamilton-Jacobi equations}.
\newblock {\em Transactions of the American Mathematical Society},
  277(1):1--42, 1983.

\bibitem{crandall1984two}
M.~G. Crandall and P.~L. Lions.
\newblock {Two approximations of solutions of Hamilton--Jacobi equations}.
\newblock {\em Mathematics of Computation}, 43(167):1--19, 1984.

\bibitem{darbon2019overcoming}
J.~Darbon, G.~P. Langlois, and T.~Meng.
\newblock Overcoming the curse of dimensionality for some Hamilton--Jacobi
  partial differential equations via neural network architectures.
\newblock {\em arXiv preprint arXiv:1910.09045}, 2019.

\bibitem{darbon2020some}
J.~Darbon and T.~Meng.
\newblock On some neural network architectures that can represent viscosity
  solutions of certain high dimensional Hamilton--Jacobi partial differential
  equations.
\newblock {\em arXiv preprint arXiv:2002.09750}, 2020.

\bibitem{darbon2016algorithms}
J.~Darbon and S.~Osher.
\newblock Algorithms for overcoming the curse of dimensionality for certain
  Hamilton--Jacobi equations arising in control theory and elsewhere.
\newblock {\em Research in the Mathematical Sciences}, 3(1):19, 2016.

\bibitem{dolgov2019tensor}
S.~Dolgov, D.~Kalise, and K.~Kunisch.
\newblock Tensor decompositions for high-dimensional Hamilton-Jacobi-Bellman
  equations.
\newblock {\em arXiv preprint arXiv:1908.01533}, 2019.

\bibitem{evans10}
L.~Evans.
\newblock {\em {Partial Differential Equations: Second Edition}}.
\newblock American Mathematical Society, 2010.

\bibitem{garcke2017suboptimal}
J.~Garcke and A.~Kr{\"o}ner.
\newblock Suboptimal feedback control of pdes by solving hjb equations on
  adaptive sparse grids.
\newblock {\em Journal of Scientific Computing}, 70(1):1--28, 2017.

\bibitem{guo2016sparse}
W.~Guo and Y.~Cheng.
\newblock A sparse grid discontinuous galerkin method for high-dimensional
  transport equations and its application to kinetic simulations.
\newblock {\em SIAM Journal on Scientific Computing}, 38(6):A3381--A3409, 2016.

\bibitem{guo2017adaptive}
W.~Guo and Y.~Cheng.
\newblock An adaptive multiresolution discontinuous {G}alerkin method for
  time-dependent transport equations in multidimensions.
\newblock {\em SIAM Journal on Scientific Computing}, 39(6):A2962--A2992, 2017.

\bibitem{guo2011local}
W.~Guo, F.~Li, and J.~Qiu.
\newblock {Local-structure-preserving discontinuous Galerkin methods with
  Lax-Wendroff type time discretizations for Hamilton-Jacobi equations}.
\newblock {\em Journal of Scientific Computing}, 47(2):239--257, 2011.

\bibitem{han2018solving}
J.~Han, A.~Jentzen, and E.~Weinan.
\newblock Solving high-dimensional partial differential equations using deep
  learning.
\newblock {\em Proceedings of the National Academy of Sciences},
  115(34):8505--8510, 2018.

\bibitem{hu1999discontinuous}
C.~Hu and C.-W. Shu.
\newblock {A discontinuous Galerkin finite element method for Hamilton--Jacobi
  equations}.
\newblock {\em SIAM Journal on Scientific Computing}, 21(2):666--690, 1999.

\bibitem{huang2019adaptive}
J.~Huang and Y.~Cheng.
\newblock An adaptive multiresolution discontinuous {G}alerkin method with
  artificial viscosity for scalar hyperbolic conservation laws in
  multidimensions.
\newblock {\em arXiv preprint arXiv:1906.00829}, 2019.

\bibitem{huang2020adaptive}
J.~Huang, Y.~Liu, W.~Guo, Z.~Tao, and Y.~Cheng.
\newblock {An adaptive multiresolution interior penalty discontinuous Galerkin
  method for wave equations in second order form}.
\newblock {\em arXiv preprint arXiv:2004.08525}, 2020.

\bibitem{huang2017quadrature}
J.~Huang and C.-W. Shu.
\newblock Error estimates to smooth solutions of semi-discrete discontinuous
  {G}alerkin methods with quadrature rules for scalar conservation laws.
\newblock {\em Numerical Methods for Partial Differential Equations},
  33(2):467--488, 2017.

\bibitem{jiang2000weighted}
G.-S. Jiang and D.~Peng.
\newblock {Weighted ENO schemes for Hamilton--Jacobi equations}.
\newblock {\em SIAM Journal on Scientific computing}, 21(6):2126--2143, 2000.

\bibitem{kang2017mitigating}
W.~Kang and L.~C. Wilcox.
\newblock Mitigating the curse of dimensionality: sparse grid characteristics
  method for optimal feedback control and hjb equations.
\newblock {\em Computational Optimization and Applications}, 68(2):289--315,
  2017.

\bibitem{ke2019alternative}
G.~Ke and W.~Guo.
\newblock {An alternative formulation of discontinous Galerkin schemes for
  solving Hamilton--Jacobi equations}.
\newblock {\em Journal of Scientific Computing}, 78(2):1023--1044, 2019.

\bibitem{kunisch2004hjb}
K.~Kunisch, S.~Volkwein, and L.~Xie.
\newblock HJB-POD-based feedback design for the optimal control of evolution
  problems.
\newblock {\em SIAM Journal on Applied Dynamical Systems}, 3(4):701--722, 2004.

\bibitem{lafon1996high}
F.~Lafon and S.~Osher.
\newblock {High order two dimensional nonoscillatory methods for solving
  Hamilton--Jacobi scalar equations}.
\newblock {\em Journal of Computational Physics}, 123(2):235--253, 1996.

\bibitem{lepsky2000analysis}
O.~Lepsky, C.~Hu, and C.-W. Shu.
\newblock {Analysis of the discontinuous Galerkin method for Hamilton--Jacobi
  equations}.
\newblock {\em Applied Numerical Mathematics}, 33(1-4):423--434, 2000.

\bibitem{li2010central}
F.~Li and S.~Yakovlev.
\newblock {A central discontinuous Galerkin method for Hamilton-Jacobi
  equations}.
\newblock {\em Journal of Scientific Computing}, 45(1-3):404--428, 2010.

\bibitem{lions1982generalized}
P.-L. Lions.
\newblock {\em Generalized solutions of Hamilton-Jacobi equations}, volume~69.
\newblock Pitman, London, 1982.

\bibitem{mallat1999wavelet}
S.~Mallat.
\newblock {\em A Wavelet Tour of Signal Processing}.
\newblock Elsevier, 1999.

\bibitem{nakamura2019adaptive}
T.~Nakamura-Zimmerer, Q.~Gong, and W.~Kang.
\newblock Adaptive deep learning for high dimensional Hamilton-Jacobi-Bellman
  equations.
\newblock {\em arXiv preprint arXiv:1907.05317}, 2019.

\bibitem{osher1988fronts}
S.~Osher and J.~Sethian.
\newblock {Fronts propagating with curvature-dependent speed: algorithms based
  on Hamilton-Jacobi formulations}.
\newblock {\em Journal of Computational Physics}, 79(1):12--49, 1988.

\bibitem{osher1991high}
S.~Osher and C.-W. Shu.
\newblock {High-order essentially nonoscillatory schemes for Hamilton--Jacobi
  equations}.
\newblock {\em SIAM Journal on Numerical Analysis}, 28(4):907--922, 1991.

\bibitem{qiu2005hermite}
J.~Qiu and C.-W. Shu.
\newblock {Hermite WENO schemes for Hamilton--Jacobi equations}.
\newblock {\em Journal of Computational Physics}, 204(1):82--99, 2005.

\bibitem{shen2010efficient}
J.~Shen and H.~Yu.
\newblock Efficient spectral sparse grid methods and applications to
  high-dimensional elliptic problems.
\newblock {\em SIAM Journal on Scientific Computing}, 32(6):3228--3250, 2010.

\bibitem{shu2007high}
C.-W. Shu.
\newblock {High order numerical methods for time dependent Hamilton-Jacobi
  equations}.
\newblock {\em Mathematics and Computation in Imaging Science and Information
  Processing, Lect. Notes Ser. Inst. Math. Sci. Natl. Univ. Singap}, 11:47--91,
  2007.

\bibitem{shu1988jcp}
C.-W. Shu and S.~Osher.
\newblock Efficient implementation of essentially non-oscillatory
  shock-capturing schemes.
\newblock {\em Journal of Computational Physics}, 77(2):439--471, 1988.

\bibitem{tao2019collocation}
Z.~Tao, Y.~Jiang, and Y.~Cheng.
\newblock An adaptive high-order piecewise polynomial based sparse grid
  collocation method with applications.
\newblock {\em arXiv preprint arXiv:1912.03982}, 2019.

\bibitem{wang2016elliptic}
Z.~Wang, Q.~Tang, W.~Guo, and Y.~Cheng.
\newblock Sparse grid discontinuous {G}alerkin methods for high-dimensional
  elliptic equations.
\newblock {\em Journal of Computational Physics}, 314:244--263, 2016.

\bibitem{yan2011local}
J.~Yan and S.~Osher.
\newblock {A local discontinuous Galerkin method for directly solving
  Hamilton--Jacobi equations}.
\newblock {\em Journal of Computational Physics}, 230(1):232--244, 2011.

\bibitem{zhang2003high}
Y.-T. Zhang and C.-W. Shu.
\newblock {High-order WENO schemes for Hamilton--Jacobi equations on triangular
  meshes}.
\newblock {\em SIAM Journal on Scientific Computing}, 24(3):1005--1030, 2003.

\end{thebibliography}

\end{document}